\documentclass{article}

\usepackage{amssymb}
\usepackage{amsmath}
\usepackage{amsthm}
\usepackage{amsfonts}
\usepackage{epsfig}
\usepackage{marvosym}

\usepackage[english]{babel}

\newtheorem{defi}{Definition}[section]
\newtheorem{prop}[defi]{Proposition}

\newtheorem{lem}[defi]{Lemma}
\newtheorem{theo}[defi]{Theorem}

\makeatletter

\@addtoreset{equation}{section}
\makeatother

\newcommand{\dx}{\mathrm d}
\newcommand{\R}{\mathbb R}

\newcommand{\E}{\mathbb E}
\newcommand{\PP}{\mathbb P}
\newcommand{\LL}{\mathbb L}

\newcommand{\N}{\mathbb N}
\newcommand{\crx}{\partial}

\def\acknow{\underline {\bf Acknowledgement}: 
We warmly thank J\'er\^ome Droniou for providing us with the article
\cite{droniou-09} and the source 
code from the corresponding deterministic
numerical scheme, and Mireille Bossy for fruitful discussions about the convergence
rate of the numerical schemes. 
}

\author{Benjamin~Jourdain$^1$, Rapha\"el~Roux$^1$}
\title{Convergence of a stochastic particle approximation for fractional
scalar conservation laws}

\begin{document}

\maketitle

{$^1$
Universit\'e Paris-Est,
CERMICS, 6 et 8 avenue Blaise Pascal, 77455 Marne-La-Vall\'ee Cedex 2, France}

\renewcommand{\thefootnote}{}
\footnotetext{\hspace{-6mm}This work was supported by the french National Research Agency (ANR)
under the program ANR-08-BLAN-0218-03 BigMC.} 
\renewcommand{\thefootnote}{\arabic{footnote}}

\abstract{In this paper, we are interested in approximating the solution to scalar
conservation laws using systems of interacting stochastic particles. The
scalar conservation law may involve a fractional Laplacian term of order
$\alpha\in (0,2]$. When $\alpha\leq 1$ as well as in the absence of this
term (inviscid case), its solution is characterized by entropic
inequalities. The probabilistic interpretation of the scalar
conservation is based on a stochastic differential equation driven by an
$\alpha$-stable process and involving a drift nonlinear in the sense of
McKean. The particle system is constructed by discretizing this equation
in time by the Euler scheme and replacing the nonlinearity by
interaction. Each particle carries a signed weight depending on its
initial position. At each discretization time we kill the couples of
particles with opposite weights and positions closer than a threshold
since the contribution of the crossings of such particles has the wrong
sign in the derivation of the entropic inequalities. We prove
convergence of the particle approximation to the solution of the
conservation law as the number of particles tends to $\infty$ whereas
the discretization step, the killing threshold and, in the inviscid
case, the  coefficient multiplying the stable increments tend to $0$ in
some precise asymptotics depending on whether $\alpha$ is larger than
the critical level $1$.}

\section*{Introduction}
We are interested in providing a numerical probabilistic
  scheme for the fractional scalar conservation law of order $\alpha$
\begin{equation}\label{eq:ldcf}
\partial_t v(t,x)+\sigma^\alpha(-\Delta)^{\frac\alpha2}v(t,x)+\partial_xA(v(t,x))=0,~~~(t,x)\in\R_+\times\R,
\end{equation}
 where $-(-\Delta)^{\frac\alpha2}$ is the fractional Laplacian operator
 of order $0<\alpha\leq2$ (defined in Section 2), and~$A$ is a function of class
$\mathcal C^1$ from $\R$ to $\R$. We also consider
 the equation obtained by letting $\sigma\to0$ in~\eqref{eq:ldci}, namely 
the inviscid conservation law
\begin{equation}\label{eq:ldci}\partial_t
  v(t,x)+\partial_xA(v(t,x))=0,~~~(t,x)\in\R_+\times\R.\end{equation}

In \cite{jourdain-02,jourdain-meleard-woyczynski-05}, these equations
are interpreted as Fokker-Planck equations associated to some
stochastic differential equations nonlinear in the sense of McKean, which can be approximated by a
particle system. We introduce an Euler time discretization of this particle
system and show the convergence of its empirical cumulative distribution
function to the solution of \eqref{eq:ldcf}.
We also study its convergence to the solution of~\eqref{eq:ldci} as the parameter $\sigma$ goes to~$0$.

Euler schemes for viscous conservation laws have already been studied in
\cite{bossy-03}, \cite{bossy-fezoui-piperno-97}, \cite{bossy-talay-96}
or \cite{bossy-talay-97}, where a convergence rate of $\frac1{\sqrt N}+\sqrt{\Delta t}$
is derived in the case $\alpha=2$, $N$ denoting the number of particles, and $\Delta t$ being
the time step.

To give the probabilistic interpretation to \eqref{eq:ldcf} we
consider the space derivative $u=\partial_xv$ of a solution~$v$
 to equation~\eqref{eq:ldcf}, which formally satisfies
\begin{equation}\label{eq:ldcf_integrée}
\partial_tu_t=-\sigma^\alpha(-\Delta)^{\frac\alpha2}u_t-\partial_x\left(A'(H*u_t)u_t\right),
\end{equation}
where $H=\mathbf1_{[0,\infty)}$ denotes the Heavyside function.
When $u_0$ is a probability measure, that is, when the initial condition
$v_0$ of Equation \eqref{eq:ldcf} is a cumulative distribution function,
Equation \eqref{eq:ldcf_integrée} 
is the Fokker-Planck equation associated to the following nonlinear
stochastic differential equation 
\begin{equation*}\begin{cases}\dx X_t&=\sigma\dx
    L_t^\alpha+A'(H*u_t(X_t))\dx t\\
u_t&=\textrm{ law of }X_t
\end{cases},\end{equation*}
where $L_t^\alpha$ is a Markov process with generator
$-(-\Delta)^{\frac\alpha2}$, namely $\sqrt2$ times a Brownian motion for $\alpha=2$, and a stable L\'evy process
with index $\alpha$ in the case $\alpha<2$, that is to say a pure jump L\'evy
process whose L\'evy measure is given by $c_\alpha\frac{\dx
  y}{|y|^{1+\alpha}}$, where $c_\alpha$ is some positive constant.

We can still give a probabilistic interpretation to Equation
\eqref{eq:ldcf} if the initial condition~$v_0$ has bounded
variation, is right continuous and not constant. Indeed, in that case $v_0$ can be written as
$v_0(x)=a+\int_{-\infty}^x\dx u_0(y)=a+H*u_0(x)$ for some finite
measure~$u_0.$ By replacing $v_0(x)$ by $\left(v_0(x)-a\right)(|u_0|(\R))^{-1}$ and
$A(x)$ by $A(a+x|u_0|(\R))(|u_0|(\R))^{-1}$ in~\eqref{eq:ldcf} ($|u_0|$
denoting the total variation of the measure $u_0,$), one
can assume without loss of generality that $a=0$ and that $|u_0|$ is a
probability measure. We denote by $\gamma=\frac{\dx u_0}{\dx|u_0|}$ the
Radon-Nikodym density of $u_0$ with respect to its total
variation. Notice that $\gamma$ takes values in~$\{\pm1\}$.

Then, Equation \eqref{eq:ldcf_integrée} is the Fokker-Planck equation associated to 
\begin{equation}\label{eq:EDS}\begin{cases}\dx X_t&=\sigma\dx
    L_t^\alpha+A'(H*\tilde P_t(X_t))\dx t\\
P&=\textrm{ law of }X
\end{cases},\end{equation}
where $\tilde P$
denotes the measure defined on the Skorokhod space $\mathcal D$ of
c\`adl\`ag functions from $[0,\infty)$ to~$\R$ by its
Radon-Nikodym density $\frac{\dx\tilde P}{\dx P}=\gamma(f(0)),$ with~$f$
the canonical process on $\mathcal D$, and $\tilde P_t$
denotes its time marginal at time $t$, {\it i.e} the measure defined by
$\tilde P_t(B)=\int_{\mathcal D}\gamma(f(0))\mathbf 1_B(f(t))\dx P(f)$,
for any~$B$ in the Borel $\sigma-$field of $\R$.

The rest of the paper is organized as follows:\\
In Section 1 we define the particle approximation for the stochastic
differential equation \eqref{eq:EDS}.\\
Section 2 is devoted to the definition of the different notions of
solutions used in the article.\\
In Section 3, we analyze the convergence of the time-discretized particle
system to the solution of the conservation law in different settings :
for both a constant or vanishing diffusion coefficient and any value of $0<\alpha\leq2$.\\
Finally, we present some numerical simulations in Section 4. Those
simulations are compared with the results of a deterministic method
described in \cite{droniou-09}.

In the following, the letter $K$ denotes some positive constant whose
value can change from line to line.

\acknow

\section{The particle approximation}\label{sect:syst_part}

In this section we construct a discretization of \eqref{eq:EDS}
consisting of both a particle approximation in order to approximate the
law of the solution and an Euler discretization to make the
particles evolve in time. The idea is to introduce $N$ particles
$X^{N,1},\hdots,X^{N,N}$ which are $N$
interacting copies of the stochastic differential equation
\eqref{eq:EDS}, where the actual law $P$ of the process is
replaced by the empirical distribution of the particles
$\frac1N\sum_{i=1}^N\delta_{X^{N,i}}$. 

In continuous time, those particles are driven by $N$ independent Brownian motions or stable
L\'evy processes with index $\alpha$ and undergo a drift given by
$A'(H*\tilde\mu_t^N(.))$, with
$\tilde\mu^N_t=\frac1N\sum_{i=1}^N\gamma(X_0^{N,i})\delta_{X^{N,i}_t}$.
The natural way to introduce the measure $\tilde\mu_t^N$ in the dynamics
is to give each particle a signed weight equal to the evaluation of $\gamma$ at the
initial position of the particle. Then, $H*\tilde\mu_t^N(x)$ is simply given
by the sum of weights of particles situated left from $x.$

The entropy solution to \eqref{eq:ldcf} has a
non-increasing total variation (see \cite{alibaud-07}), which can be
interpreted probabilistically as a compensation of merging sample paths having opposite
signs. For a more precise statement in the case $\alpha=2$, see Lemma~2.1
in~\cite{jourdain-02}. It is thus natural to adapt this behavior in our
particle approximation by killing
any merging couple of particles with opposite signs.

In \cite{jourdain-02} Jourdain proves, for $\alpha=2$ in continuous
time, the convergence of the particle system to the solution of the
nonlinear stochastic differential equation through a
propagation-of-chaos result. Moreover, the convergence of the signed cumulative
distribution function~$H*\tilde\mu_t^N$ to the solution to
Equation~\eqref{eq:ldcf} is also proved, as well as convergence to the
solution to the
inviscid equation as $\sigma\to0$.
 In~\cite{jourdain-meleard-woyczynski-05} the same results are generalized
to the case $1<\alpha<2,$ assuming $\gamma=1$ in the case of a vanishig viscosity. However, to our knowledge there is even no
existence result for the particle system in continuous time when
$\alpha\leq1$, since the driving L\'evy process is somehow weaker than the
nonlinear drift.

In discrete time, the probability of
seeing two particles actually merging is 0. To adapt the murders from
the continuous time setting, we thus kill, at each time step, any couple
of particles with opposite signs separated by a distance smaller
than a given threshold~$\varepsilon_N$ going to zero as $N$ goes to
$\infty$. Though, one has to be careful, since one can
have more than two particles lying in a small interval of length~$\varepsilon_N.$ 
Precisely, the particles are killed in the following way: kill the
leftmost couple of particles at consecutive positions separated
by a distance smaller than the threshold $\varepsilon_N$ and with
opposite signs. Then, recursively apply the same algorithm to the
remaining particles. This can be done with a computational cost of order $\mathcal O(N)$.
The essential properties satisfied by this killing procedure are the following:
\begin{itemize}
\item to each killed particle is attached another killed
  particle, which has opposite signs and lies at a distance at most
  $\varepsilon_N$ of the first particle.
\item after the killing there is no couple of particles with
 opposite signs in a distance smaller than~$\varepsilon_N.$
\item the exchangeability of the particles is preserved.
\item after the murder, the quantity $H*\tilde\mu_t^{N}(X_t^{N,i})$
  remains the same for any surviving particle.
\end{itemize}

We are going to describe the killed processes by a couple $(f,\kappa)$
in the space $\mathcal K=\mathcal D\times[0,\infty]$ 
of c\`adl\`ag functions $f$ from
$[0,\infty)$ to $\R$ endowed with a death time
$\kappa\in[0,\infty]$. The space $\mathcal K$ is endowed with the
product metric
$d((f,\kappa_f),(g,\kappa_g))=d_S(f,g)+|\arctan(\kappa_f)-\arctan(\kappa_g)|,$
where $d_S$ is the Skorokhod metric on $\mathcal D,$ so that $(\mathcal K,d)$ is a
complete metric space.
It could seem more natural to consider the space $\mathcal
D([0,\infty),\R\cup\{\crx\})$ of paths taking values in $\R$ endowed with a
cemetery point $\crx.$ However the corresponding topology is too strong to
prove Proposition \ref{prop:tension}.

The precise description of the process is the following: each particle
will be represented by a couple $(X^{N,i},\kappa_i^N)\in\mathcal K$. Let
$(X_0^i)_{i\in\N}$ be a sequence of independent random variables with
common distribution~$|u_0|$ and let $h_N>0$ denote the time step of the Euler scheme.
At time 0, kill the particles according to the preceding rules,
that is to say, set $\kappa_i^N=0$ for killed particles, which will not
move anymore. Those particles
will not be taken into account anymore.
Now, by induction, suppose that the particle system has been defined up to time~$kh_N,$ and
kill the particles according to the preceding rules ({\it i.e.} set
$\kappa_i^N=kh_N$ and $X_t^{N,i}=X_{kh_N}^{N,i}$ for all $t\geq kh_N$, if the particle
with index $i$ is one of those). Then let
the particles still alive evolve up to time $(k+1)h_N$ according to
 $$\dx X_t^{N,i}=A'\left(\frac1N\sum_{\kappa_j^N>kh_N}\gamma(X_0^j)\mathbf1_{X_{kh_N}^{N,j}\leq X_{kh_N}^{N,i}}\right)\dx
t+\sigma_N\dx L_t^i,$$ where $(L^i)_{i\in\N}$ is a sequence of independent
$\alpha$-stable L\'evy processes for $\alpha<2$, or a sequence of
independent copies of $\sqrt2$ times Brownian motion, which are independent of the sequence~$(X_0^i)_{i\in\N}$.
The particle system is thus well-defined, by induction.

Let $\mu^N=\frac1N\sum_{i=1}^N\delta_{(X^{N,i},\kappa_i^N)}\in\mathcal P(\mathcal K)$ be the empirical distribution of the
particles. For a probability measure $Q$ on~$\mathcal K$ and $t\geq0$, we define a
signed measure $\tilde Q_t$ on $\R$ by:
$$\tilde Q_t(B)=\int_{\mathcal K}\mathbf1_B(f(t))\mathbf1_{\kappa>t}\gamma(f(0))\dx Q(f,\kappa),$$
for any $B$ in the Borel $\sigma-$field of $\R.$
With these notations, on the interval $[kh_N,(k+1)h_N),$ a particle,
provided it is still alive, satisfies
$$\dx X_t^{N,i}=A'\left(H*\tilde\mu_{kh_N}^N\left(X_{kh_N}^{N,i}\right)\right)\dx t
+\sigma_N\dx L_t^i.$$
Notice that the sum of the weights of alive particles $\tilde\mu_t^N(\R)=\frac1N\sum_{\kappa_i^N>t}\gamma(X_0^i)$ is constant
in time, since the particles are killed by couples of opposite signs.

\section{Notion of solutions}

In this section, we recall the different notions of solutions that are
associated to the equations~\eqref{eq:ldcf} and~\eqref{eq:ldci}. Indeed,
due to the shock-creating term $\partial_x(A(u_t))$, the notion of weak solution is too weak,
and does not provide uniqueness when the diffusion term is not
regularizing enough. The best suited notion in
those cases is the notion of
entropy solution.

In \cite{kruzhkov-70}, Kruzhkov shows that for $v_0\in\LL^\infty((0,\infty))$ existence and uniqueness hold for
entropy solutions to \eqref{eq:ldci}, defined as functions
$v\in\LL^\infty((0,\infty)\times\R)$ satisfying, for any smooth convex
function $\eta,$ any nonnegative smooth function $g$ with compact support
on $[0,\infty)\times\R$ and any $\psi$
satisfying $\psi'=\eta'A'$, the entropic inequality
\begin{equation}
\label{eq:entropie}
\int_\R
\eta(v_0)g_0+\int_0^\infty\left(\int_\R\eta(v_t)\partial_tg_t+\psi(v_t)\partial_xg_t\right)\dx
t\geq0.
\end{equation}
It is well known that this entropy solution can be obtained as the limit of weak solutions to
\eqref{eq:ldcf} as $\sigma\to0$ in the case $\alpha=2$. 

Weak solutions to~\eqref{eq:ldcf} (see
\cite{jourdain-02}) are defined as functions
$v\in\LL^\infty((0,\infty)\times\R)$ satisfying, for all smooth
functions $g$ with compact support in~$[0,\infty)\times\R$,
\begin{equation}\label{eq:faible}
\int_\R v_0g_0+\int_0^\infty\int_\R v_t\partial_tg_t\dx
t-\sigma^\alpha\int_0^\infty\int_\R v_t(-\Delta)^{\frac\alpha2}
g_t\dx t+\int_0^\infty\int_\R A(v_t)\partial_xg_t\dx
t=0.
\end{equation}

For $\alpha<2,$ we denote by $(-\Delta)^{\frac\alpha2}$ the fractional symmetric differential operator of order
$\alpha,$ that can be defined through the Fourier transform:
$$\widehat{(-\Delta)^{\frac\alpha2} u}(\xi)=|\xi|^\alpha\hat u(\xi).$$ 
An equivalent definition for $(-\Delta)^{\frac\alpha2}$ uses an
integral representation
$$(-\Delta)^{\frac\alpha2}u(x)=c_\alpha\int_\R\frac{u(x+y)-u(x)-\mathbf
  1_{|y|\leq r}u'(x)y}{|y|^{1+\alpha}}\dx y$$
for any $r\in(0,\infty)$ and some fixed constant $c_\alpha$ (see \cite{droniou-imbert-06}),
depending on the definition of the Fourier transform.

It has
been proven in \cite{jourdain-02} and \cite{jourdain-meleard-woyczynski-05} that existence and
uniqueness holds for weak solutions of~\eqref{eq:ldcf}, for
$1<\alpha\leq2$.
However, for $0<\alpha\leq1$, the diffusive term of order $\alpha$ in \eqref{eq:ldcf} is
somehow dominated by the shock-creating term, which is of order 1, so that a weak formulation
does not ensure uniqueness for the solution. We thus have to strengthen
the notion of solution, and use entropy
solutions to \eqref{eq:ldcf}, defined in~\cite{alibaud-07}
as functions $v$ 
in $\LL^\infty((0,\infty)\times\R)$ 
satisfying the relation 
\begin{align}\label{eq:entropifr}
&\int_0^\infty\eta(v_0)g_0+\int_0^\infty\int_\R\left(\eta(v_t)\partial_tg_t+\psi_t(v_t)\partial_xg_t\right)\dx
t\nonumber\\
&\quad+c_\alpha\int_0^\infty\int_\R\int_{\{|y|>
  r\}}\eta'(v_t(x))\frac{v_t(x+\sigma y)-v_t(x)}{|y|^{1+\alpha}}g_t(x)\dx
y\dx x\dx t\\
&\quad+c_\alpha\int_0^\infty\int_\R\int_{\{|y|\leq
  r\}}\eta(v_t(x))\frac{g_t(x+\sigma y)-g_t(x)-\sigma y\partial_xg_t(x)}{|y|^{1+\alpha}}\dx
y\dx x\dx t\geq0\nonumber
\end{align}
for any $r>0$, any nonnegative smooth
function $g$ with compact support in $[0,\infty)\times\R$, any smooth
convex function $\eta:\R\rightarrow\R$ and any $\psi$ satisfying
$\psi'=\eta'A'$. Notice that from the convexity of $\eta$, the entropic
formulation~\eqref{eq:entropifr}
 for a parameter $r$ implies the entropic
formulation with parameter $r'>r$. Also notice, using the functions
$\eta(x)=\pm x$ that an entropy solution to \eqref{eq:ldcf} is a weak
solution to~\eqref{eq:ldcf}.

In \cite{alibaud-07}, Alibaud shows that existence and uniqueness hold
for entropy solutions of \eqref{eq:ldcf} provided that the initial
condition $v_0$ lies in $\LL^\infty(\R)$. The entropy solution then lies
in the space $\mathcal C([0,\infty),\LL^1(\frac{\dx x}{1+x^2}))$. He also proves that the
entropy solution to~\eqref{eq:ldcf}
converges  to the entropy solution
to \eqref{eq:ldci} in the space $\mathcal
C([0,T],\LL^1_{\textrm{loc}}(\R))$ as~$\sigma\to0$.

\section{Statement of the results}

The aim of this article is to prove the three following convergence
result, each one corresponding to a particular setting.
\begin{theo}\label{theo:conv_alpha<1}
Assume $0<\alpha\leq1.$ Let $\sigma_N\equiv\sigma$ be a constant sequence.
Let $\varepsilon_N$ and $h_N$ be two vanishing sequences satisfying the inequalities
$$N^{-\lambda}\leq4\sup_{[-1,1]}|A'|h_N\leq\varepsilon_N,\textrm{ and
}N^{-1/\alpha}\leq N^{-1/\lambda}\varepsilon_N$$for
some positive $\lambda.$ For $\alpha=1$,
 also assume $h_N\leq\varepsilon_N N^{-1/\lambda}$. It holds for any $T>0$,
$$\lim_{N\to\infty}\int_0^T\E\left\|H*\tilde \mu^N_t-v_t\right\|_{\LL^1\left(\frac{\dx
      x}{1+x^2}\right)}\dx t=0,$$
where $v_t$ denotes the entropy solution to
 the fractional conservation law \eqref{eq:ldcf}. 
\end{theo}
\begin{theo}\label{theo:conv_sigma0}
Let $\varepsilon_N,$ $h_N$ and $\sigma_N$ be three vanishing
  sequences such that
  $$N^{-\lambda}\leq4\sup_{[-1,1]}|A'|h_N\leq\varepsilon_N$$for some
  $\lambda>0$. If
  $\alpha>1$, also assume
  $\sigma_N\leq\varepsilon_N^{1-\frac1\alpha}N^{-\frac1\lambda}$. Then, for any $T>0,$
$$\lim_{N\to\infty}\int_0^T\E\left\|H*\tilde \mu^N_t-v_t\right\|_{\LL^1\left(\frac{\dx
      x}{1+x^2}\right)}\dx t=0,$$
where $v_t$ denotes the entropy solution to
 the inviscid conservation law \eqref{eq:ldci}.
\end{theo}
The additional assumption for $\alpha>1$ comes from the fact that in
this case, the dominant term is the diffusion, while in the limit there
is no diffusion anymore. The assumption ensures that the diffusion is
weak enough not to perturb the approximation. For $\alpha\leq1,$ the
dominant term is the drift, as in the limit, so that no additional
condition is needed.
\begin{theo}\label{theo:conv_alpha>1}
Assume $1<\alpha\leq2$. Let $\sigma_N\equiv\sigma$ be a constant
sequence, and let $\varepsilon_N$ and $h_N$ be two vanishing
sequences.
It holds for any $T>0$,
$$\lim_{N\to\infty}\int_0^T\E\left\|H*\tilde \mu^N_t-v_t\right\|_{\LL^1\left(\frac{\dx
      x}{1+x^2}\right)}\dx t=0,$$
where $v_t$ denotes the weak solution to the fractional conservation law
\eqref{eq:ldcf}.
\end{theo}

In order to prove those three theorems, we will have to control the probability of
seeing particles merging. In the case $\alpha<2$, this is mainly due to the
conjunction of the small jumps of the stable process and the drift
coefficient, while the large jumps of the stable term do not play an essential role.
 As a consequence, for $\alpha<2$, we consider another family of
evolutions coinciding with the Euler scheme on the time
discretization grid, for which we
consider differently the jumps which are smaller or larger than a given
threshold~$r$. The
choice of this parameter has to be linked to the parameter $r$
appearing in the entropic formulation \eqref{eq:entropifr}, since they
play a similar role: the third term in
\eqref{eq:entropifr} corresponds to the effect of jumps larger than $r$
in the driving L\'evy process and the fourth
term corresponds to jumps smaller than $r$. This evolution is
designed so that on the first half of each time step, the process will
evolve according to the drift and the small jumps, and on the second half
of each time step, it will evolve according to the large jumps. More precisely, let $$\nu^i(\dx y,\dx
t)=\sum_{\Delta L^i_t\neq0}\delta_{(\Delta L_t^i,t)}$$ be the jump measure associated
to the L\'evy process $L^i$ and let 
$$\tilde\nu^i(\dx y,\dx t)=\nu^i(\dx y,\dx t)-c_\alpha\frac{\dx y\dx t}{|y|^{1+\alpha}}$$ be the corresponding
compensated measure, so that
$$L_t^i=\int_{(0,t]\times\{|y|>r\}}y\nu^i(\dx y,\dx
t)+\int_{(0,t]\times\{|y|\leq r\}}y\tilde\nu^i(\dx y,\dx t),$$ where the
right hand side does not depend on $r$.
We define the process $X^{N,i,r}$ by
$$X^{N,i,r}=X_0^i+\sigma_NL^{N,i,r}+\sigma_N\Lambda^{N,i,r}+\mathcal
A^{N,i},$$ where
\begin{itemize}
\item $L_t^{N,i,r}$ is the large jumps part defined by 
$$L^{N,i,r}_t=\int_{(0,a(t)]\times\{|y|>r\}}y\nu^i(\dx
 y,\dx s ),$$
where $a(t)=\begin{cases}kh_N&\textrm{ for }t\in[kh_N,(k+1/2)h_N]\\
kh_N+2(t-(k+1/2)h_N)&\textrm{ for
}t\in[(k+1/2)h_N,(k+1)h_N]\end{cases}.$ This process is constant on
intervals $[kh_N,(k+1/2)h_N]$ and behaves like a L\'evy process with jump
measure $\mathbf1_{|y|>r}\frac{2c_\alpha\dx y}{|y|^{1+\alpha}}$ on intervals $[(k+1/2)h_N,(k+1)h_N]$.
\item $\Lambda^{N,i,r}_t$ is the small jumps part, defined by 
$$\Lambda^{N,i,r}_t=\int_{(0,b(t)]\times\{|y|\leq r\}}\tilde\nu^i(\dx y,\dx s),$$ 
where $b(t)=\begin{cases}kh_N+2(t-kh_N)&\textrm{ for }t\in[kh_N,(k+1/2)h_N]\\(k+1)h_N
&\textrm{ for }t\in[(k+1/2)h_N,(k+1)h_N]\end{cases}$. This term behaves like a L\'evy process with jump
measure $\mathbf1_{|y|\leq r}\frac{2c_\alpha\dx y}{|y|^{1+\alpha}}$ on
intervals $[kh_N,(k+1/2)h_N]$ and is constant on
intervals $[(k+1/2)h_N,(k+1)h_N]$. Notice that the process
$\Lambda^{N,i,r}$ is a martingale.
\item $\mathcal A^{N,i}$ is the drift part, which satisfies
  $\mathcal A_0^{N,i}=0$, is
  constant over each interval 
$[(k+1/2)h_N,(k+1)h_N]$, and
evolves as a piecewise
affine process with derivative
$2A'(H*\tilde\mu_{kh_N}^N(X_{kh_N}^{N,i}))$ on intervals
$[kh_N,(k+1/2)h_N].$
\end{itemize}
One can check that for any $r$, the process
$(X^{N,1,r},\hdots,X^{N,N,r})$ is equal to $(X^{N,1},\hdots,X^{N,N})$ on
the time discretization grid up to killing time.
Conditionally on the positions of the particles at time~$kh_N,$ the
particles evolve independently on $[kh_N,(k+1)h_N]$, and the evolution
on $[kh_N,(k+1/2)h_N]$ is independent of the evolution on $[(k+1/2)h_N,(k+1)h_N].$
Since the entropic formulation~\eqref{eq:entropifr} with parameter $r$
is stronger than the one with parameter $r'\geq r,$ we have to make the
parameter $r$ tend to zero in order to prove the entropic formulation
for any parameter. However, this convergence has
to satisfy some conditions with respect to $N,$ $h_N$ and
$\varepsilon_N$. We will explain later why a suitable sequence $r_N$
exists under the conditions given
in the statement of Theorem~\ref{theo:conv_alpha<1}.

In order to prove Theorems \ref{theo:conv_alpha<1} and
\ref{theo:conv_sigma0}, we introduce $\mu^{N,r}$ the empirical distribution of the processes
$(X^{N,i,r},\kappa_i^N)$:
$$\mu^{N,r}=\frac1N\sum_{i=1}^N\delta_{(X^{N,i,r},\kappa_i^N)}\in\mathcal P(\mathcal
K),$$ and by $\pi^{N,r}$ the law of~$\mu^{N,r}$.

The following proposition is the first step in the proof of Theorems
\ref{theo:conv_alpha<1}, \ref{theo:conv_sigma0}
and \ref{theo:conv_alpha>1}.
\begin{prop}\label{prop:tension}
\begin{itemize}
\item Assume $\alpha<2$.
For any bounded
 sequences $(h_N)$, $(\sigma_N)$ and $(\varepsilon_N)$, and for any
 sequence $(r_N)$,
 the family of probability measures $(\pi^{N,r_N})_{N\in\N}$ is tight in
$\mathcal P(\mathcal P(\mathcal K))$.
\item Denote by $\pi^N$ the law of $\mu^N.$
For any bounded
 sequences $(h_N)$, $(\sigma_N)$ and $(\varepsilon_N)$,
 the family of probability measures $(\pi^N)_{N\in\N}$ is tight in
$\mathcal P(\mathcal P(\mathcal K))$.
\end{itemize}
\end{prop}
\begin{proof}
We first check the tightness of the family $(\pi^{N,r_N})_{N\in\N}$.

As stated in \cite{sznitman-89}, checking the tightness of the sequence
$\pi^{N,r_N}$ boils down to checking the tightness of the sequence
$(\textrm{Law}(X^{N,1,r_N},\kappa_1^N)).$
Owing to the product-space structure, we can
check tightness for~$X^{N,1,r_N}$ and  $\kappa_1^N$ separately.

Of course, tightness for $\kappa_1^N$ is straightforward since it lies
on the compact space $[0,\infty]$, and it is
enough to check tightness for the laws of the path $(X^{N,1,r_N})$.
For simplicity, we will assume that $A=0$, which is not restrictive
since $A'$ is a bounded function so that the perturbation induced by~$A$
belongs to a compact subset of the space of continuous functions, from Ascoli's theorem (also
notice that the addition functional from $\mathcal
D\times\mathcal C([0,\infty))$ to $\mathcal D$ is
continuous).
We use Aldous' criterion to prove tightness (see \cite{aldous-78}).
First, the sequences $(X_0^{N,1,r_N})_{N\in\N}$ and $(\sup_{[0,T]}|\Delta
X^{N,1,r_N}|)_{N\in\N}$ are tight, since $(X_0^{N,1,r_N})$ is constant
in law and $\left(\sup_{[0,T]}|\Delta
X^{N,1,r_N}|\right)_N$ is dominated by the identically distributed
sequence $$\left(\left(\sup_N\sigma_N\right)\sup_{[0,T+\sup_Nh_N]}\left|\Delta
L^1\right|\right)_N.$$

Then let $\tau_N$ be a stopping time of the natural filtration of
$X^{N,1,r_N}$ taking finitely many values, and let~$(\delta_N)_{N\in\N}$ be a sequence of positive
numbers going to $0$ as $N\to\infty.$ One can write
\begin{align}\label{eq:continuite}
\PP\left(\left|X_{\tau_N+\delta_N}^{N,1,r_N}-X_{\tau_N}^{N,1,r_N}\right|\geq\varepsilon\right)
\leq&\PP\left(\sigma_N\left|\Lambda_{\tau_N+\delta_N}^{N,1,r_N}
-\Lambda_{\tau_N}^{N,1,r_N}\right|\geq\varepsilon/2\right)\nonumber\\
&+
\PP\left(\sigma_N\left|L_{\tau_N+\delta_N}^{N,1,r_N}-L_{\tau_N}^{N,1,r_N}\right|\geq\varepsilon/2\right)\\
\leq&\PP\left(\sup_{t\in[0,\delta_N]}\sigma_N|\mathcal L_t^{\leq
    r_N}|\geq\varepsilon/2\right)+
\mathbb
P\left(\sup_{t\in[0,\delta_N]}\sigma_N|\mathcal
  L_t^{>r_N}|\geq\varepsilon/2\right),\nonumber
\end{align}
where $$\mathcal L^{\leq r}_t=\int_{(0,t]\times\{|y|\leq r\}}y\tilde\nu(\dx y,\dx
  t)\textrm{~~and~~}\mathcal L^{>r}_t=\int_{(0,t]\times\{|y|>r\}}y\nu(\dx y,\dx
  t),$$ the measure $\nu$ being the jump measure of some L\'evy
process $\mathcal L$ with L\'evy measure $\frac{2c_\alpha\dx
  y}{|y|^{1+\alpha}}$, and $\tilde\nu$ is the compensated measure of $\nu$.
Now, using the maximal inequality for the
martingale~$(\mathcal L_t^{\leq r_N})_{t\in[0,\delta_N]},$
 noticing that $(\mathcal L_{\delta_N}^{\leq r})_{r\in[0,1]}$ is also a martingale, we deduce
\begin{align*}
\PP\left(\sup_{t\in[0,\delta_N]}|\mathcal L_t^{\leq
    r_N}|\geq\varepsilon/2\sigma_N\right)&\leq\sup_{r\in[0,\sup_Nr_N]}\PP\left(\sup_{t\in[0,\delta_N]}|\mathcal L_t^{\leq
    r}|\geq\varepsilon/2\sigma_N\right)\\
&\leq2\sigma_N\varepsilon^{-1}\sup_{r\in[0,\sup_Nr_N]}\E\left(|\mathcal
  L_{\delta_N}^{\leq r}|\right)\\
&=2\sigma_N\varepsilon^{-1}\E\left(|\mathcal L_{\delta_N}^{\leq\sup_Nr_N}|\right)\\
&\underset{N\to\infty}\longrightarrow0.
\end{align*}
For the large jumps parts, one writes,
\begin{align*}
\PP\left(\sup_{t\in[0,\delta_N]}|\mathcal
  L_t^{>r_N}|\geq\varepsilon/2\sigma_N\right)&\leq\PP\left(\sup_{t\in[0,\delta_N]}|\mathcal
  L_t|+\sup_{t\in[0,\delta_N]}|\mathcal
  L_t^{\leq r_N}|\geq\varepsilon/2\sigma_N\right)\\
&\leq\PP\left(\sup_{t\in[0,\delta_N]}|\mathcal
  L_t|\geq\varepsilon/4\sigma_N\right)+\PP\left(\sup_{t\in[0,\delta_N]}|\mathcal
  L_t^{\leq r_N}|\geq\varepsilon/4\sigma_N\right)\\
&\underset{N\to\infty}\to0.
\end{align*}
As a consequence, the family $(\textrm{Law}(X^{N,1,r_N}))_{N\in\N}$ is tight in $\mathcal D$.

Thus, the family $(\pi^{N,r_N})_{N\in\N}$ is tight.

The proof is essentially the same for the tightness of
$(\pi^N)_{N\in\N}$, with a few simplifications, since we do not treat
separately large and small jumps. It also adapts in the case $\alpha=2$,
since the Gaussian distribution has thinner tails than the
$\alpha-$stable distribution for $\alpha<2.$
\end{proof}

The use of the path space $\mathcal K$ instead of
$\mathcal D([0,\infty),\R\cup\{\crx\})$ for a cemetery point $\crx$ is crucial in the proof of Proposition
\ref{prop:tension}, since in the latter case, we need to control the
jumps occuring close to the death time in order to prove tightness. The following
example is illustrative: if we consider a sequence $f_n$ of paths starting at $0$, jumping to 1 at time
$1-1/n$, and being killed at time $1$, then~$f_n$ does not converge in
$\mathcal D([0,\infty),\R\cup\{\crx\})$, while it does in $\mathcal
K$.

The following lemma deals with the initial condition of the particle system.

\begin{lem}\label{lem:cond_init}
If $\pi^\infty$ is the limit of some subsequence of $\pi^N$ or $\pi^{N,r_N}$, then for
$\pi^\infty$-almost all $Q$, for all $A$ in the Borel $\sigma-$field of $\R$,
\begin{equation}\label{eq:Q0}
Q_0(A):=\int_\R\mathbf1_{\kappa>0}\mathbf1_{f(0)\in A}\dx
Q(f,\kappa)=|u_0|(A).
\end{equation}
\end{lem}
In particular, $\kappa$ is $Q-$almost surely positive for
$\pi^\infty$-almost all $Q$.
\begin{proof}
In a first time, we control the probability of seeing a partincle dying
within a short time.

Let us write the Hahn decomposition $u_0^+-u_0^-$ of the measure
$u_0$, the measures $u_0^+$ and $u_0^-$ being positive measures supported by two disjoint
sets $B^+$ and $B^-.$ From the inner regularity of the measure $u_0^+,$
for any $\delta>0$, one can
find a closed set $F^+\subset B^+$ such that $u_0^+(F^+)\geq
u_0^+(B^+)-\delta.$ The complement set $O^-=(F^+)^c$ is then an open
subset of $\R$, which can thus be decomposed as a countable union of disjoint
open intervals $O^-=\bigcup_{m=1}^\infty]a_m,b_m[$. For a large
enough integer $M$, and for $\varepsilon_\delta>0$
small enough, the set
$O^\delta=\bigcup_{m=1}^M]a_m+\varepsilon_\delta,b_m-\varepsilon_\delta[$
is such that $u_0^-(O^\delta)\geq u_0^-(O^-)-\delta.$
Consequently, we can write $\R$ as a partition $$\R=F^+\cup(B^-\cap
O^\delta)\cup\mathcal B^\delta,$$ where $\mathcal B^\delta=(F^+\cup(B^-\cap
O^\delta))^c$ has
small measure $|u_0|(\mathcal B^\delta)\leq
2\delta,$ particles starting in $F^+$ have a positive sign, and
particles starting in $(B^-\cap O^\delta)$ have a negative sign.
Let $N$ be large enough to ensure $\varepsilon_N\leq\varepsilon_\delta/3.$
The distance between any element of~$F^+$ and any element
of $O^\delta$ is larger than~$\varepsilon_\delta$. As a
consequence, if the particles with index~$i$ and $j$ kill
each other before time $\tau$, then either one of them started in~$\mathcal B^\delta$, or one of the particles $i$ and $j$
moved by a distance larger than $\varepsilon_\delta/3$.
This writes 
\begin{align*}
\sharp\left\{i,\kappa_i^N<\tau_\delta\right\}&=2\sharp\left\{(i,j),i<j,X^{N,i,r_N}\textrm{
  and }X^{N,j,r_N}\textrm{ kill
each other}\right\}\\
&\leq2\sharp\left\{i,X_0^{N,i,r_N}\in\mathcal
B^\delta\textrm{ or }\sup_{t\in[0,\tau]}|X_t^{N,i,r_N}-X_0^i|\geq\varepsilon_\delta/3\right\}.
\end{align*}
As a consequence,
if~$\tau_\delta>0$ is small enough so that
$\PP(\sup_{t\in[0,\tau_\delta]}|X_t^{N,i,r_N}-X_0^{N,i,r_N}|\geq\varepsilon_\delta/3)\leq\delta$
(this can be achieved using an adaptation of~\eqref{eq:continuite}),
it holds
\begin{align*}
\PP(\kappa_1^N<\tau_\delta)=\frac1N\E\left(\sharp\left\{i,\kappa_i^N<\tau_\delta\right\}\right)&\leq\frac2N\E\left(\sharp\left\{i,X_0^{N,i,r_N}\in\mathcal
B^\delta\textrm{ or }\sup_{t\in[0,\tau_\delta]}|X_t^{N,i,r_N}-X_0^i|\geq\varepsilon_\delta/3\right\}\right)
\\&\leq2\PP\left(X_0^i\in\mathcal B^\delta\right)+2\PP\left(\sup_{t\in[0,\tau_\delta]}|X_t^{N,i,r_N}-X_0^{N,i,r_N}|\geq\varepsilon_\delta/3\right)\\
&\leq6\delta.
\end{align*}
Consequently,
$$\E^{\pi^\infty}(Q(\kappa<\tau))\leq\underset
N{\lim\inf}~\E^{\pi^N}(Q(\kappa<\tau))=\underset N{\lim\inf}~\PP(\kappa_1^N<\tau)\underset{\tau\to0}\to0.$$
Thus for $\pi^\infty$-almost all $Q$, $\kappa$ is $Q-$almost surely positive.
As a consequence, for any bounded continuous function $\varphi$,
\begin{align*}
\E^{\pi^\infty}\left|\int_{\mathcal
    K}\mathbf1_{\kappa>0}\varphi(f(0))\dx Q(\kappa,f)-\int_\R\varphi\dx |u_0|\right|
&=\E^{\pi^\infty}\left|\int_{\mathcal K}\varphi(f(0))\dx
  Q(\kappa,f)-\int_\R\varphi\dx |u_0|\right|\\
&=\lim_N\E^{\pi^N}\left|\int_{\mathcal K}\varphi(f(0))\dx
  Q(\kappa,f)-\int_\R\varphi\dx |u_0|\right|=0,
\end{align*}from the law of large numbers.

\end{proof}

The main step in the proof of Theorems~\ref{theo:conv_alpha<1},~\ref{theo:conv_sigma0} and~\ref{theo:conv_alpha>1}  is the following
proposition:

\begin{prop}\label{prop:convergence}
Let $\varepsilon_N$ and $h_N$ be vanishing sequences.
\begin{itemize}
\item If $\sigma_N$ is a constant sequence and $0<\alpha\leq1$, suppose
  $N^{-1/\alpha}\leq N^{-1/\lambda}\varepsilon_N$ and $N^{-\lambda}\leq4\sup_{[-1,1]}|A'|h_N\leq\varepsilon_N$
 for some
positive $\lambda$. If
  $\alpha=1,$ also assume $h_N\leq N^{-1/\lambda}\varepsilon_N.$ Then, there
exists a sequence $(r_N)$ of positive real numbers, such that the limit of any converging subsequence
 of $\pi^{N,r_N}$ gives full measure
to the set$$\{Q\in\mathcal P(\mathcal K), H*\tilde Q_t(x) \textrm{ is the entropy solution to \eqref{eq:ldcf}}\}.$$
\item
Let $\sigma_N$ be a vanishing sequence and assume $N^{-\lambda}\leq4\sup_{[-1,1]}|A'|h_N\leq\varepsilon_N$
 for some
positive~$\lambda$. If $1<\alpha\leq2$, also assume
  $\sigma_N\leq\varepsilon_N^{1-\frac1\alpha}N^{-1/\lambda}$. Then
$$\{Q\in\mathcal P(\mathcal K), H*\tilde Q_t(x) \textrm{ is
  the entropy solution to \eqref{eq:ldci}}\}$$ is given full measure by
any limit of a converging subsequence of $\pi^{N,r_N},$ for a well chosen
sequence~$(r_N),$ in the case $\alpha<2$, and by any limit of a
converging subsequence of $\pi^N$ if $\alpha=2.$
\item If $\sigma_N$ is a constant sequence and $1<\alpha\leq2$, the limit of any converging subsequence
 of $\pi^N$ gives full measure
to the set$$\{Q\in\mathcal P(\mathcal K), H*\tilde Q_t(x) \textrm{ is
  the weak solution to \eqref{eq:ldcf}}\}.$$
\end{itemize}
\end{prop}

Proposition \ref{prop:convergence} will be proved in Section
\ref{sect:preuve}. We first admit it to end the proofs
of Theorems~\ref{theo:conv_alpha<1},~\ref{theo:conv_sigma0}  and~\ref{theo:conv_alpha>1}. We need the following lemma.
\begin{lem}\label{lem:approx_r_N}Let $\alpha<2$ and $r_N$ be a sequence of positive numbers going to
  zero. Then it holds, for any $T>0,$
$$\underset{N\to\infty}
\lim\int_0^T\E\|H*\tilde\mu_t^N-H*\tilde\mu_t^{N,r_N}\|_{\LL^1\left(\frac{\dx
      x}{1+x^2}\right)}\dx t=0.$$
\end{lem}
\begin{proof}
It holds, by exchangeability of the particles,
\begin{align*}
\int_0^T\E\|H*\tilde\mu_t^N-H*\tilde\mu_t^{N,r_N}\|_{\LL^1\left(\frac{\dx
      x}{1+x^2}\right)}\dx t&\leq\E\int_0^T\int_\R\frac1N\sum_{\kappa_i^N>t}\left|\mathbf1_{X_t^{N,i}\leq x}-\mathbf1_{X_t^{N,i,r_N}\leq
  x}\right|\frac{\dx x\dx t}{x^2+1}\\
&\leq\int_0^T\E\left(\mathbf1_{\kappa_1^N\geq
    t}\left|X_t^{N,1}-X_t^{N,1,r_N}\right|\wedge\pi\right)\dx t.
\end{align*}
This last quantity goes to zero, since the processes $X^{N,1}$ and
$X^{N,1,r_N}$ coincide on the discretization grid, whose mesh
vanishes. Indeed, for $t\in[kh_N,(k+1)h_N)$
\begin{align}\label{eq:continuite2}
\E\left(\mathbf1_{\kappa_1^N>t}|X_t^{N,1,r_N}-X_t^{N,1}|\wedge\pi\right)
&\leq\E\left(\mathbf1_{\kappa_1^N>t}|X_t^{N,1,r_N}-X_{kh_N}^{N,1}|\wedge\pi\right)
+\E\left(\mathbf1_{\kappa_1^N>t}|X_t^{N,1}-X_{kh_N}^{N,r_N}|\wedge\pi\right)\nonumber\\
&\leq
Kh_N^{1/2}.
\end{align}
For this last estimate, we used, for an $\alpha-$stable L\'evy process $L$, the inequality
$$\E\left(|L_t|\wedge1\right)\leq K\E\left(|L_t|^{\alpha/2}\right)=Kt^{1/2}.$$
\end{proof}
From Lemma \ref{lem:approx_r_N}, it is sufficient to show
$\underset{N\to\infty}\lim\int_0^T\E\|H*\tilde\mu_t^{N,r_N}-v_t\|_{\LL^1(\frac{\dx
    x}{1+x^2})}\dx t$
in order to prove Theorems~\ref{theo:conv_alpha<1} and \ref{theo:conv_sigma0}.
\begin{proof}[Proof of Theorems
  \ref{theo:conv_alpha<1}-\ref{theo:conv_sigma0}-\ref{theo:conv_alpha>1}]
We write the proof for Theorems \ref{theo:conv_alpha<1} and
\ref{theo:conv_sigma0} in the case $\alpha<2$.
The proof of Theorem \ref{theo:conv_sigma0} with $\alpha=2$ and Theorem \ref{theo:conv_alpha>1} is the same, with $\pi^N$
replacing $\pi^{N,r_N}$.

Let $\gamma^k$ be a Lipschitz continuous approximations of $\gamma$, with
$\PP(\gamma(X_0^1)\neq \gamma^k(X_0^1))\leq1/k$ (see \cite{jourdain-02},
Lemma 2.5, for a construction of such a $\gamma^k$). We have, by exchangeability of the particles,
\begin{align}\label{eq:approx_r_N}
&\E\int_0^T\int_\R\left|H*\tilde\mu_t^{N,r_N}(x)-v_t(x)\right|\frac{\dx
  x\dx t}{x^2+1}\nonumber\\
\leq&\E\int_0^T\int_\R\mathbf1_{\kappa_1^N>
    t}H(x-X_t^{N,1,r_N})\left|\gamma(X_0^{N,1,r_N})-\gamma^k(X_0^{N,1,r_N})\right|\frac{\dx x\dx t}{x^2+1}\\
&+\E^{\pi^N}\left(\int_0^\infty\int_\R\left|\int_\mathcal
    K\mathbf1_{\kappa>t}H(x-f(t))\gamma^k(f(0))\dx Q(f,\kappa)-v_t(x)\right|\frac{\dx x}{x^2+1}\right).\nonumber
\end{align}
From the assumption on $\gamma^k,$
the first term in the right hand side of \eqref{eq:approx_r_N} is
smaller than
$2\pi/k$ which vanishes as $k$ goes to $\infty.$
The bounded function$$Q\mapsto\int_0^T\int_\R\left|\int_{\mathcal
    K}\mathbf1_{\kappa>t}H(x-f(t))\gamma^k(f(0))\dx
  Q(f,\kappa)-v_t(x)\right|\frac{\dx x\dx t}{x^2+1}$$is continuous.
From Proposition \ref{prop:convergence}, the second term in the right
hand side of \eqref{eq:approx_r_N} converges, as~$N$ goes to $\infty$ to
$$\E^{\pi^\infty}\left(\int_0^T\int_\R\left|\int_{\mathcal K}\mathbf1_{\kappa>t}H(x-f(t))\left(\gamma^k(f(0))
-\gamma(f(0))\right)\dx Q(f,\kappa)\right|\frac{\dx
x}{x^2+1}\right).$$
This terms goes to zero as $k$ tend to infinity using the argument of
the begining of the proof with~$X^{N,1,r_N}$ replaced by the canonical
process $y$.

\end{proof}

\subsection{Proof of Proposition \ref{prop:convergence}}\label{sect:preuve}
 This section is devoted to the proof of Proposition
 \ref{prop:convergence}. Since the hardest part of this proof is the
 first two items, we do not give all details for the third item and for
 the second one in the case $\alpha=2$. Indeed, for these two last
 settings, the separation of small jumps and large jump is not necessary
 for the proof.

Let $r_N$ be a sequence of positive real numbers, going to zero as $N\to\infty$, which will be
explicited later. Let $r>0$ and $c$ be reals numbers, $\eta$ a smooth convex
function, $\psi$ a primitive of $A'\eta'$ and~$g$ a smooth compactly
supported nonnegative
function. We define the function 
$\varphi_t(x)=\int_{-\infty}^xg_t(y)\dx y$. Note that $\varphi$ is
smooth, and nondecreasing with respect to the space variable. We consider a subsequence of $\pi^{N,r_N}$, still
denoted~$\pi^{N,r_N}$ for simplicity, which converges to a limit
$\pi^\infty$. 
We want to prove that, for $\pi^\infty-$almost all $Q$, the function
$H*\tilde Q_t$ satisfies the entropy formulation associated to the
corresponding case.

One can write, for any $k\geq0$ and $t\in]kh_N,(k+1)h_N]$
\begin{align*}&\PP\left(\exists i,j,\kappa_i^N\wedge\kappa_j^N>t,X_t^{N,i,r}
=X_t^{N,j,r}\right)=\E\left(\PP\left(\exists
  i,j,\kappa_i^N\wedge\kappa_j^N>t,X_t^{N,i,r}=X_t^{N,j,r}\bigg|(X_{kh_N}^{N,q})_q\right)\right)\\
=&\E\left(\PP\left(\exists i,j,\kappa_i^N\wedge\kappa_j^N>t,
\sigma_NZ^{i,j,k,N}_t=X_{kh_N}^{N,j}-X_{kh_N}^{N,i}+\mathcal
A^{N,j}_t-\mathcal A^{N,i}_t\bigg|(X_{kh_N}^{N,q})_q\right)\right),
\end{align*}
where we denote 
$$Z_t^{i,j,N,k}=\Lambda_t^{N,i,r}-\Lambda_{kh_N}^{N,i,r}-\Lambda_t^{N,j,r}
+\Lambda_{kh_N}^{N,j,r}+L_t^{N,i,r}-L_{kh_N}^{N,i,r}-L_t^{N,j,r}+L_{kh_N}^{N,j,r}.$$
From the conditional independence of the processes $L^{N,i,r}$,
$L^{N,j,r}$, $\Lambda^{N,i,r}$ and $\Lambda^{N,j,r}$, the random variable
$Z_t^{i,j,N,k}$ has a density. As a consequence, since the process $\mathcal
A^{N,j}_t-\mathcal A^{N,i}_t$ is deterministic on $[kh_N,(k+1)h_N]$ conditionally to $(X_{kh_N}^{N,q})_q$, the
above probability is zero, meaning that for all time $t>0$, the
alive particles $X_t^{N,i,r_N}$ almost surely have distinct positions. As a consequence, the function
$\eta\left(H*\tilde\mu_t^{N,r_N}(x)\right)$
is the cumulative distribution function of the signed
measure $$\xi_t^N=\sum_{\kappa_i^N>t}w_t^i\delta_{X_t^{N,i,r_N}},$$
where 
\begin{align*}
w_t^i&=\mathbf1_{\kappa_i^N>t}\left(\eta\left(\frac1N\sum_{\substack{\kappa_j^N>t\\X_t^{N,j,r_N}\leq
      X_t^{N,i,r_N}}}\gamma(X_0^j)\right)
-\eta\left(\frac1N\sum_{\substack{\kappa_j^N>t\\X_t^{N,j,r_N}<X_t^{N,i,r_N}}}\gamma(X_0^j)\right)\right)\\
&=\mathbf1_{\kappa_i^N>t}\left(\eta\left(H*\tilde\mu_t^{N,r_N}\left(X_t^{N,i,r_N}\right)\right)
-\eta\left(H*\tilde\mu_t^{N,r_N}\left(X_t^{N,i,r_N}-\right)\right)\right).
\end{align*}
Let $(\zeta_m)_{m\in\N}$ be the
increasing sequence of times which are either a jump time for some
$L^{N,i,r_N}$ ({\it
  i.e.} a jump of size $>r_N$ for $X^{N,i,r_N}$) or either a time of the form~$kh_N/2$.
One has
\begin{align}\label{eq:egalite}
-\left<\xi_0^N,\varphi_0\right>=&\sum_{m=1}^\infty\left<\xi_{\zeta_m}^N,\varphi_{\zeta_m}\right>
-\left<\xi_{\zeta_{m-1}}^N,\varphi_{\zeta_{m-1}}\right>\nonumber\\
=&\sum_{\kappa_i^N>0}\sum_{m=1}^\infty
w_{\zeta_{m-1}}^i\left(\varphi_{\zeta_m}\left(X_{\zeta_m-}^{N,i,r_N}
  \right)
-\varphi_{\zeta_{m-1}}\left(X_{\zeta_{m-1}}^{N,i,r_N}\right)\right)\\
&+\sum_{\kappa_i^N>0}\sum_{m=1}^\infty\left(w_{\zeta_m}^i\varphi_{\zeta_m}\left(X_{\zeta_m}^{N,i,r_N}
  \right)-w_{\zeta_{m-1}}^i\varphi_{\zeta_m}\left(X_{\zeta_m-}^{N,i,r_N}\right)\right).\nonumber
\end{align}
Notice that these infinite sums are actually finite, since the function
$\varphi_t$ is identically zero when $t$ is large enough, and since the
process $(L^{N,1,r_N},\hdots,L^{N,N,r_N})$ has a finite number of jumps on bounded intervals.

We consider the first term in the right hand side of \eqref{eq:egalite}.
Denote by $\displaystyle\nu^{i,r}=\sum_{\Delta
X_t^{N,i,r}\neq0}\delta_{(\Delta L_t^{N,i,r}+\Delta\Lambda_t^{N,i,r},t)}$ the jump measure associated to
$L^{N,i,r}+\Lambda^{N,i,r}$, and by 
$$\tilde \nu^{i,r}(\dx y,\dx t)=\nu^{i,r}(\dx y,\dx
t)-2c_\alpha\left(\chi_t^N\mathbf1_{|y|\leq r}+(1-\chi_t^N)\mathbf1_{|y|> r}\right)\frac{\dx
  y\dx t}{|y|^{1+\alpha}}$$
 its compensated measure, where
$\chi^N_t=\sum_{k=0}^\infty\mathbf 1_{[kh_N,(k+1/2)h_N)}(t)$. Let us
apply It\=o's Formula on the interval 
$(\zeta_{m-1},\zeta_m)$. If
$\zeta_{m-1}=kh_N$ for some integer $k$, then $\zeta_m=(k+1/2)h_N$, and
almost surely $X_{(k+\frac12)h_N-}^{N,i,r}=X_{(k+\frac12)h_N}^{N,i,r}$
holds. As a consequence
\begin{align*}&\varphi_{(k+\frac12)h_N}\left(X_{(k+\frac12)h_N-}^{N,i,r}\right)-\varphi_{kh_N}\left(X_{kh_N}^{N,i,r}\right)
\\=&\int_{kh_N}^{(k+\frac12)h_N}\partial_t\varphi_t(X_t^{N,i,r})\dx
t+2\int_{kh_N}^{(k+\frac12)h_N}\partial_x\varphi_t(X_t^{N,i,r})
A'\left(H*\tilde\mu_{kh_N}^{N,r}(X_{kh_N}^{N,i,r})\right)\dx t.\\
&+\int_{(kh_N,(k+1/2)h_N)}\int_{\{|y|\leq r\}}\left(\varphi_t(X_{t-}^{N,i,r}+\sigma_Ny)
-\varphi_t(X_{t-}^{N,i,r})-\sigma_Ny\partial_x\varphi_t(X_{t-}^{N,i,r})\right)\nu^{i,r}(\dx y,\dx t)\\
&+\sigma_N\int_{(kh_N,(k+1/2)h_N)}\partial_x\varphi_t(X_{t-}^{N,i,r})\left(\int_{\{|y|\leq
  r\}}y\tilde \nu^{i,r}(\dx y,\dx t)\right).
\end{align*}
If $\zeta_{m-1}$ is not of the form $kh_N$, then the process $X^{N,i,r}$ is constant on the
interval $[\zeta_{m-1},\zeta_m)$, and one has
$\varphi_{\zeta_m}(X_{\zeta_m-}^{N,i,r})-\varphi_{\zeta_{m-1}}(X_{\zeta_{m-1}}^{N,i,r})
=\int_{\zeta_{m-1}}^{\zeta_m}\partial_t\varphi_t(X_t^{N,i,r})\dx t$.
Summing over all the intervals $(\zeta_{m-1},\zeta_m),$ Equation
\eqref{eq:egalite} writes,
denoting $\tau_t=\max\{\zeta_m,\zeta_m\leq t\}$,

\begin{align}\label{eq:ineg}
-\left<\xi_0^N,\varphi_0\right>
&=\sum_{\kappa_i^N>0}\int_0^\infty w_{\tau_t}^i\left(\partial_t\varphi_t(X_t^{N,i,r_N})
+2\chi_t^N\partial_x\varphi_t(X_t^{N,i,r_N})A'\left(H*\tilde\mu_{\tau_t}^{N,r_N}(X_{\tau_t}^{N,i,r_N})\right)\right)\dx t\nonumber\\
&\quad+c_\alpha\sum_{\kappa_i^N>0}\int_0^\infty w_{\tau_t}^i\chi_t^N\int_{\{|y|\leq
  r_N\}}\left(\varphi_t(X_t^{N,i,r_N}+\sigma_Ny)
-\varphi_t(X_t^{N,i,r_N})
-\sigma_Ny\partial_x\varphi_t(X_t^{N,i,r_N})\right)\frac{2\dx y\dx
t}{|y|^{1+\alpha}}\nonumber\\
&\quad+\sum_{\kappa_i^N>0}~\sum_{\substack{\textrm{large jump}\\\textrm{at }\zeta_m}}w_{\zeta_m}^i\varphi_{\zeta_m}(X_{\zeta_m}^{N,i,r_N})
-w_{\zeta_{m-1}}^i\varphi_{\zeta_m}(X_{\zeta_m-}^{N,i,r_N})
\nonumber\\
&\quad+\sum_{\kappa_i^N>0}~\sum_{\substack{\zeta_m\textrm{ of the}\\
\textrm{form }kh_N}}(w_{\zeta_m}^i-w_{\zeta_{m-1}}^i)\varphi_{\zeta_m}(X_{\zeta_m}^{N,i,r_N})\\
&\quad+\sum_{\kappa_i^N>0}\sum_{\substack{\zeta_m\textrm{ of the}\\\textrm{form }(k+1/2)h_N}}
(w_{\zeta_m}^i-w_{\zeta_{m-1}}^i)\varphi_{\zeta_m}(X_{\zeta_m}^{N,i,r_N})\nonumber\\
&\quad+M_N.\nonumber
\end{align}
Here, the third, fourth and fifth terms correspond to the second term
in the right hand side of \eqref{eq:egalite}, and $M_N$ is a martingale term given by
$$M_N=\sum_{\kappa_i^N>0}\int_0^\infty w_{\tau_t}^i\chi_t^N\int_{\{|y|\leq r_N\}}
\left(\varphi_t(X_{t-}^{N,i,r_N}+\sigma_Ny)-\varphi_t(X_{t-}^{N,i,r_N})\right)\tilde\nu^{i,r_N}(\dx y,\dx t).$$
Equation \eqref{eq:ineg} can be rewritten
$$T_N^1=T_N^2+T_N^3+T_N^4+T_N^5+M_N,$$
where $T_N^1=-\left<\xi_0^N,\varphi_0\right>$, $T_N^2$ is the sum of the two first terms in the right-hand-side of
\eqref{eq:ineg}, $T_N^3$ is the third one, $T_N^4$ the fourth one and $T_N^5$
the fifth one.

The four following Lemmas, whose proofs are postponed to Section~\ref{sect:preuves_lemmes}
 deal with the asymptotic behavior of the terms $M_N$, $T_N^2-T_N^1$, $T_N^3$
and $T_N^4$.
\begin{lem}\label{lem:martingale}
It holds $$\E|M_N|^2\leq\frac{K\sigma_N^2r_N^{2-\alpha}}{N}$$for some positive constant $K$.
The equivalent term in the case $\alpha=2$,
$$M_N=\sigma_N\sum_{\kappa_i^N>0}\int_0^\infty w_{\tau_t}^i\partial_x\varphi(X_t^{N,i})\dx L_t^i,$$satifies the same estimate :
$$\E|M_N|^2\leq K\frac{\sigma_N^2}N.$$
\end{lem}
\begin{lem}\label{lem:T2}
\begin{itemize}
\item It holds
\begin{align*}
\E\bigg|&-T_N^1+T_N^2+\int_0^\infty\int_\R\left(\eta(H*\tilde\mu_t^{N,r_N})\partial_tg_t
+2\chi_t^N\psi(H*\tilde\mu_t^{N,r_N})\partial_xg_t\right)\dx
t+\int_\R g_0\eta(H*\tilde\mu_0^{N,r_N})\dx x\\
&\left.+2c_\alpha\int_0^\infty\chi_t^N\int_\R\int_{\{|y|\leq
  r_N\}}\eta(H*\tilde\mu_t^{N,r_N}(x))\left(g_t(x+\sigma_Ny)-g_t(x)-\sigma_Ny\partial_xg_t(x)\right)\frac{\dx y\dx
  x\dx t}{|y|^{1+\alpha}}\right|\underset{N\to\infty}\to0.
\end{align*}
\item If $r_N\leq1/\sigma_N$, then $$\left|2c_\alpha\int_0^\infty\chi_t^N\int_\R\int_{\{|y|\leq
  r_N\}}\eta(H*\tilde\mu_t^{N,r_N}(x))\left(g_t(x+\sigma_Ny)-g_t(x)-\sigma_Ny\partial_xg_t(x)\right)\frac{\dx y\dx
  x\dx t}{|y|^{1+\alpha}}\right|\leq K\sigma_N^\alpha.$$
\end{itemize}
\end{lem}

The following lemma gives two estimates for the term $T_N^3$, the first
being useful for a constant viscosity $\sigma_N\equiv\sigma$, and the
second for vanishing viscosity $\sigma_N\to0$.
\begin{lem}\label{lem:T3}
\begin{itemize}
\item The error term
$$\E\left|T_N^3+2c_\alpha\int_0^\infty(1-\chi_t^N)\int_\R\int_{\{|y|>r_N\}}
\!\!\!\!\!\!\eta'(H*\tilde\mu^{N,r_N}_t(x))\left(H*\tilde\mu_t^{N,r_N}(x+\sigma_Ny)-
H*\tilde\mu_t^{N,r_N}(x)\right)g_t(x)\frac{\dx y\dx x\dx
t}{|y|^{1+\alpha}}\right|$$
vanishes if $N^{-1}r_N^{-\alpha}$ goes to~0.
\item  It holds $$\E|T_N^3|\leq K(\sigma_Nr_N^{1-\alpha}+\sigma_N^\alpha).$$
\end{itemize}
\end{lem}

\begin{lem}\label{lem:T4}
One has $\E|T_N^4|\underset{N\to\infty}\to0.$
\end{lem}

We now have to control the probability for the last remaining term $T_N^5$ to
be negative.
If there is no crossing of particles with opposite signs between $kh_N$
and $(k+1/2)h_N$, for any $k$, then $T_N^5\geq~\!\!0$.
Indeed, let $X_{(k+1/2)h_N}^{N,i_1,r_N}\leq\hdots\leq X^{N,i_q,r_N}_{(k+1/2)h_N}$ be a maximal sequence of consecutive
particles with same sign. The sequence
$\left(\varphi_{(k+1/2)h_N}(X_{(k+1/2)h_N}^{N,i_l,r_N})\right)_{l=1,\hdots,q}$ is thus a nondecreasing sequence,
and from the convexity of $\eta$ and the fact that no particles with
opposite signs cross, $(w_{(k+1/2)h_N}^{i_l})_{l=1,\hdots,q}$ is the nondecreasing
reordering of $(w_{kh_N}^{i_l})_{l=1,\hdots,q}$.
Thus, from Lemma~\ref{lem:max_somme} below, 
$\sum_{\kappa_i^N>kh_N}(w_{(k+1/2)h_N}^i-w_{kh_N}^i)\varphi_{(k+1/2)h_N}(X_{(k+1/2)h_N}^{N,i,r_N})$ is nonnegative.
 It is thus sufficient to control the probability that two particles
 with opposite signs cross between $kh_N$ and $(k+1/2)h_N$. Since
 after the murder there is no couple of particles with opposite signs
 separated by a smaller distance than $\varepsilon_N,$ this does not happen as soon as no
particle drift by more than $\varepsilon_N/4$ and no particle is moved
by more than $\varepsilon_N/4$ by the small jumps. The drift on half a time step is smaller
than $\sup_{[-1,1]}|A'|h_N$ which is assumed to be smaller than $\varepsilon_N/4$.
We control the contribution of the small jumps in the following lemma:

\begin{lem}\label{lem:controle_sauts}
Let $B_N$ be the event $$B_N=\left\{\forall k\leq T/h_N,\forall
i,\sigma_N\left|\Lambda_{(k+1/2)h_N}^{i,r_N}-\Lambda_{kh_N}^{i,r_N}\right|
\leq\varepsilon_N/4\right\},$$ so that no
crossing of particles with opposite signs between
$kh_N$ and $(k+1/2)h_N$ occurs on $B_N$. One  has, for $\alpha<2$,
$$\PP(B_N)\geq\left(1-e^{Kh_Nr_N^{-\alpha}-\varepsilon_N/4\sigma_Nr_N}\right)^{NT/h_N},$$

For $\alpha=2,$ we define the event $B_N$ by $$B_N=\left\{\forall k\leq T/h_N,\forall
i,\sigma_N\left|L_{(k+1)h_N}^i-L_{kh_N}^i\right|
\leq\varepsilon_N/4\right\}.$$It holds
$$\PP(B_N)\geq\left(1-Ke^{-\varepsilon_N^2/(32h_N\sigma_N^2)}\right)^{NT/h_N}.$$
\end{lem}
The proof will be given in Section \ref{sect:preuves_lemmes}.

We now gather all the previous information to prove that, depending on
the considered case, the entropic
formulation or the weak formulation holds almost surely.

\begin{enumerate}
\item  {\it Constant viscosity $\sigma_N\equiv\sigma$, with index $0<\alpha\leq1$.}

Define, for $Q\in\mathcal P(\mathcal K)$,
\begin{align*}
F_N^r(Q)=&\int_\R\eta(H*\tilde Q_0)g_0+\int_0^\infty
\int_\R\left(\eta(H*\tilde Q_t)\partial_tg+2\chi_t^N\psi(H*\tilde Q_t)\partial_xg\right)\dx
t\\
&\quad+2c_\alpha\int_0^\infty(1-\chi_t^N)\int_\R\int_{\{|y|>
  r\}}\eta'(H*\tilde Q_t(x))(H*\tilde Q_t(x+\sigma_N y)-H*\tilde Q_t(x))g_t(x)\frac{\dx y\dx x\dx t}{|y|^{1+\alpha}}\\
&\quad+2c_\alpha\int_0^\infty\chi_t^N\int_\R\int_{\{|y|\leq
  r\}}\eta(H*\tilde Q_t(x))(g_t(x+\sigma_N y)-g_t(x)-\sigma_N y\partial_xg_t(x))\frac{\dx
y\dx x\dx t}{|y|^{1+\alpha}}
\end{align*}
and
\begin{align*}
F^r(Q)=&\int_\R\eta(H*\tilde Q_0)g_0+\int_0^\infty\int_\R\left(\eta(H*\tilde Q_t)\partial_tg
+\psi(H*\tilde Q_t)\partial_xg\right)\dx t\\
&\quad+c_\alpha\int_0^\infty\int_\R\int_{\{|y|>
  r\}}\eta'(H*\tilde Q_t(x))(H*\tilde Q_t(x+\sigma y)-H*\tilde Q_t(x))g_t(x)\frac{\dx y\dx x\dx t}{|y|^{1+\alpha}}\\
&\quad+c_\alpha\int_0^\infty\int_\R\int_{\{|y|\leq
  r\}}\eta(H*\tilde Q_t(x))(g_t(x+\sigma y)-g_t(x)-\sigma y\partial_xg_t(x))\frac{\dx
y\dx x\dx t}{|y|^{1+\alpha}}.
\end{align*}
Notice that from the convexity of $\eta,$ one has 
$$\eta'(H*\tilde Q_t(x))(H*\tilde Q_t(x+\sigma y)-H*\tilde Q_t(x))\leq
\eta(H*\tilde Q_t(x+\sigma y))-\eta(H*\tilde Q_t(x)),$$
 so that for any $0<r\leq r'$, it holds $F^r\leq F^{r'}$ and $F_N^r\leq F_N^{r'}$.

From Equation \eqref{eq:ineg}, it holds, for $N$ large enough so that $r_N\leq r$,
$$F^r_N(\mu^{N,r_N})\geq
F^{r_N}_N(\mu^{N,r_N})=T_N^5+\left(-T_N^1+T_N^2+T_N^3+T_N^4+M_N+F^{r_N}_N(\mu^{N,r_N})\right).$$
From the assumptions made on $\varepsilon_N$ and $h_N$ one can construct a
 sequence $r_N$ such that $N^{-1/\alpha}=o(r_N)$,
 $h_Nr^{-\alpha}_N=o\left(\varepsilon_Nr_N^{-1}\right)$ and $\frac
 N{h_N}e^{-\varepsilon_N/4\sigma r_N}\to0$.
Indeed, set $r_N=\varepsilon_NN^{-1/2\lambda}$. Then it holds
$N^{-1/\alpha}r_N^{-1}\leq KN^{-1/2\lambda}$ and
$\frac{h_N}{\varepsilon_N}r_N^{1-\alpha}=
h_N\varepsilon_N^{-\alpha}N^{(\alpha-1)/2\lambda}$, which vanishes for
any value of $\alpha$.
Then $\frac N{h_N}$ goes
to infinity at the rate of a power of $N$, and $\varepsilon_N/r_N=N^{1/2\lambda}$ as
well. Thus, $\frac N{h_N}e^{-\varepsilon_N/4\sigma r_N}$ vanishes.

 As a consequence, from Lemmas \ref{lem:martingale}, \ref{lem:T2}, \ref{lem:T3} and
 \ref{lem:T4}, $\E\left|-T_N^1+T_N^2+T_N^3+T_N^4+M_N+F^{r_N}_N(\mu^{N,r_N})\right|$ vanishes
as $N$ tends to infinity, and the event $B_N$ defined in Lemma \ref{lem:controle_sauts} is such that $\PP(B_N)\to1$.
On the event $B_N$, $T_N^5$ is almost-surely
nonnegative, so that, from the uniform boundedness of~$F^r_N$ with
respect to $N$,
$\E^{\pi^{N,r_N}}(F^r_N(Q)^-)=\E(F^r_N(\mu^{N,r_N})^-)$
 goes to~$0$. To show that the entropic formulation holds almost
 surely, we need a continuous approximation of $F_N^r$ and $F^r.$ We define $F^{r,\delta}$ and
$F_N^{r,\delta}$ by replacing every occurrence of $H*\tilde Q_t$
in the definitions of $F^r$ and $F_N^r$ by
$\int_{\mathcal K}\mathbf1_{\kappa>t}H(.-f(t))\gamma^\delta(f(0))\dx
Q(f,\kappa)$, 
where~$\gamma^\delta$ is a Lipschitz continuous approximation of $\gamma$, with
$\PP(\gamma(X_0^1)\neq \gamma^\delta(X_0^1))\leq\delta$ (see \cite{jourdain-02},
Lemma~2.5, for the construction of $\gamma^\delta$).
Then, for any fixed $\delta$ and $r$,
the family $\{F^{r,\delta}\}\cup\{F_N^{r,\delta},N\in\N\}$ is equicontinuous for the
topology of weak convergence.
Indeed, let~$Q^k$ be a sequence of probability measures on $\mathcal K$ converging to $Q$
as $k$ goes to infinity. From the continuity of the application
$\mathcal K\to\R, (f,\kappa)\mapsto\mathbf1_{\kappa>0}f(0)$, $Q_0^k$ converges
weakly to $Q_0$ (where $Q_0$ and $Q^k_0$ are defined as in~\eqref{eq:Q0}), 
and from the continuity of the applications $\mathcal K\to\R,(f,\kappa)\mapsto
\mathbf1_{\kappa>t}\gamma^\delta(f(0))\mathbf1_{f(t)\leq y}$ on the set $\{(f,\kappa)\in\mathcal K,
f(t)=f(t-),f(t)\neq y\}$, for all $t$ in the complement of the countable
set $\{t\in[0,\infty),Q(\{f(t)\neq f(t-)\}\cup\{\kappa=t\})>0\},$
the quantity $\int_{\mathcal K}\mathbf1_{\kappa>t}H(.-f(t))\gamma^\delta(f(0))\dx Q^k(f,\kappa)$
converges almost everywhere to $\int_{\mathcal K}\mathbf1_{\kappa>t}H(.-f(t))\gamma^\delta(f(0))\dx Q(f,\kappa)$.
From Lebesgue's bounded convergence theorem, we deduce that
$$\sup_N|F_N^{r,\delta}(Q^k)-F_N^{r,\delta}(Q)|+|F^{r,\delta}(Q^k)-F^{r,\delta}(Q)|\underset{k\to\infty}\to0$$
yielding equicontinuity for $\{F^{r,\delta}\}\cup\{F_N^{r,\delta},N\in\N\}$.
Moreover, since the sequence $\chi_t^N$ converges~$*$-weakly to $1/2$ in the space
$\LL^\infty((0,\infty))$, $F_N^{r,\delta}$ converges pointwise to
$F^{r,\delta}$ as $N$ goes to infinity. Ascoli's
theorem thus implies that $F_N^{r,\delta}$
converges uniformly on compact sets to $F^{r,\delta}.$ 
From the weak convergence of $\pi^{N,r_N}$ to~$\pi^\infty,$ one thus
deduces $$\E^{\pi^{N,r_N}}[F_N^{r,\delta}(Q)^-]\underset{N\to\infty}\to\E^{\pi^\infty}[F^{r,\delta}(Q)^-].$$
Moreover, for any $t>0$, any $y$, and any probability measure $Q$
satisfying $Q_0=|u_0|$ (with $Q_0$ defined as in \eqref{eq:Q0}), which holds true for $\pi^\infty-$almost all $Q$
from Lemma \ref{lem:cond_init}, it holds 
$$\left|H*\tilde Q_t(y)-\int_{\mathcal
  K}\mathbf1_{\kappa>t}H(y-f(t))\gamma^\delta(f(0))\dx
Q(f,\kappa)\right|\leq\int_\R|\gamma-\gamma^\delta|\dx|u_0|\leq\delta,$$
yielding convergence to $0$ for
$\E^{\pi^\infty}|F^r(Q)^--F^{r,\delta}(Q)^-|+\E^{\pi^{N,r_N}}|F^r_N(Q)^--F^{r,\delta}_N(Q)^-|$
as $\delta$ goes to~$0$, uniformly in $N$.
 As a consequence, writing 
\begin{align*}
\E^{\pi^\infty}(F^r(Q)^-)\leq&\E^{\pi^\infty}|F^r(Q)^--F^{r,\delta}(Q)^-|
+\left|\E^{\pi^\infty}(F^{r,\delta}(Q)^-)-\E^{\pi^{N,r_N}}(F_N^{r,\delta}(Q)^-)\right|\\
&+\E^{\pi^{N,r_N}}|F^{r,\delta}_N(Q)^--F^r_N(Q)^-|+\E^{\pi^{N,r_N}}(F_N^r(Q)^-)
\end{align*}
we deduce that $F^r(Q)$ is nonnegative for~$\pi^\infty-$almost all
$Q$. We just have to notice that Lemma~\ref{lem:cond_init} yields that, $\pi^\infty-$almost surely, $H*\tilde Q_0=v_0$ to conclude
that the entropy formulation holds~$\pi^\infty-$almost surely.

\item {\it Vanishing viscosity $\sigma_N\to0$.}

We define
$$F_N^r(Q)=\int_\R\eta(H*\tilde Q_0)g_0+\int_0^\infty
\int_\R\left(\eta(H*\tilde Q_t)\partial_tg+2\chi_t^N\psi(H*\tilde
  Q_t)\partial_xg\right)\dx t$$
and
$$F(Q)=\int_\R\eta(H*\tilde Q_0)g_0+\int_0^\infty\int_\R\left(\eta(H*\tilde Q_t)\partial_tg
+\psi(H*\tilde Q_t)\partial_xg\right)\dx t.$$
Regularized versions $F_N^{r,\delta}$ and $F^\delta$ of $F_N^r$ and
$F$ are also considered using the function $\gamma^\delta$ instead of~$\gamma$. 
In the case $\alpha<2,$ the same arguments as above, using the second parts of Lemmas~\ref{lem:T2} and~\ref{lem:T3} will show that the
entropy formulation holds $\pi^\infty-$almost surely for $H*\tilde
Q_t$, provided
there exists a sequence $r_N$ such that
$\frac{\sigma_N^2r_N^{2-\alpha}}N$ and $\sigma_Nr_N^{1-\alpha}$ vanish, $r_N\leq\sigma_N^{-1},$
 $h_Nr^{-\alpha}_N=o(\varepsilon_N(\sigma_Nr_N)^{-1})$ and $\frac
 N{h_N}e^{-\varepsilon_N/4\sigma_N r_N}\to0$.
\begin{itemize}
\item For $\alpha\leq1,$ any sequence $r_N$ vanishing at a very quick
  rate will fit.

\item For $\alpha>1$, since we assumed $\sigma_N\leq
\varepsilon_N^{1-\frac1\alpha}N^{-1/\lambda}$ these conditions are
satisfied by the sequence
 $r_N=\frac{\varepsilon_N}{\sigma_N}N^{-\frac{\alpha}{2\lambda(\alpha-1)}}.$
\end{itemize}
In the case $\alpha=2$, It\=o's formula writes
\begin{align*}\varphi_{(k+1)h_N}\left(X_{(k+1)h_N}^{N,i}\right)-\varphi_{kh_N}\left(X_{kh_N}^{N,i}\right)
=&\int_{kh_N}^{(k+1)h_N}\partial_t\varphi_t(X_t^{N,i})\dx
t\\
&+2\int_{kh_N}^{(k+1)h_N}\partial_x\varphi_t(X_t^{N,i})
A'\left(H*\tilde\mu_{kh_N}^N(X_{kh_N}^{N,i})\right)\dx t.\\
&+\sigma_N^2\int_{(kh_N,(k+1)h_N)}\int_{\{|y|\leq
  r\}}\partial_x^2\varphi_t(X_t^{N,i})\dx t\\
&+\sigma_N\int_{(kh_N,(k+1)h_N)}\partial_x\varphi_t(X_t^{N,i})\dx L_t^i.
\end{align*}
The three first terms are treated as in the case $\alpha<2$, and the
stochastic integral is dealt with using Lemma \ref{lem:martingale}.
For the entropic inequality to holds, we need to control the crossing of
particles with opposite sign. From Lemma \ref{lem:controle_sauts}, if $\frac
N{h_N}e^{-\frac{\varepsilon_N^2}{32\sigma_N^2 h_N}}$ goes to zero, then no
crossing occurs. Since our assumptions yield
$h_N\sigma_N^2\leq\varepsilon_N^2N^{-1/\lambda}$ and $N/h_N\leq
KN^{1+\lambda}$, this condition holds true.

\item {\it Constant viscosity $\sigma_N\equiv\sigma,$ with index $1<\alpha\leq2$.}

In this case, since we want to derive a weak formulation, we do not need to consider separately large and small
jumps. As a consequence it is enough to study the process $X_t^{N,i}.$

Let $g$ be a smooth function with compact support, and define for $Q\in\mathcal P(\mathcal K)$,
$$F(Q)=\int_\R H*\tilde
Q_0g_0+\int_0^\infty\int_\R H*\tilde Q_t\partial_tg_t\dx
t-\sigma^\alpha\int_0^\infty\int_\R H*\tilde
Q_t(-\Delta)^{\frac\alpha2}g_t\dx t+\int_0^\infty\int_\R A(H*\tilde Q_t)\partial_xg_t.$$

Let $\varphi_t(x)=\int_{-\infty}^xg_t(y)\dx y$.
One has
\begin{align*}
-\frac1N\sum_{\kappa_i^N>0}\gamma(X_0^{N,i})\varphi_0(X_0^{N,i})
=&-\frac1N\sum_{k=0}^\infty\sum_{\kappa_i^N=(k+1)h_N}\gamma(X_0^{N,i})\varphi_{(k+1)h_N}(X_{(k+1)h_N}^{N,i})\\
&+\frac1N\sum_{k=0}^\infty\sum_{\kappa_i^N>kh_N}
\gamma(X_0^{N,i})\left(\varphi_{(k+1)h_N}(X_{(k+1)h_N}^{N,i})-\varphi_{kh_N}(X_{kh_N}^{N,i})\right).
\end{align*}
From It\=o's formula, in the case $\alpha<2,$ when $\kappa_i^N>kh_N$, 
\begin{align}\label{eq:Ito}
&\varphi_{(k+1)h_N}(X_{(k+1)h_N}^{N,i})-\varphi_{kh_N}(X_{kh_N}^{N,i})\\
=&\int_{kh_N}^{(k+1)h_N}\partial_t\varphi_t(X_t^{N,i})\dx
t+\int_{kh_N}^{(k+1)h_N}\partial_x\varphi_t(X_t^{N,i})A'\left(H*\tilde\mu_{kh_N}(X_{kh_N}^{N,i})\right)\dx
t\nonumber\\
&+c_\alpha\int_{(kh_N,(k+1)h_N)}\int_\R\left(\varphi_t(X_t^{N,i}+\sigma y)
-\varphi_t(X_t^{N,i})-\mathbf1_{\{|y|\leq r\}}\sigma y\partial_x\varphi_t(X_t^{N,i})\right)\frac{\dx y\dx t}{|y|^{1+\alpha}}\nonumber\\
&+\int_{(kh_N,(k+1)h_N)}\int_\R\left(\varphi_t(X_{t-}^{N,i}+\sigma y)
-\varphi_t(X_{t-}^{N,i})\right)\tilde\nu^i(\dx y,\dx t).\nonumber
\end{align}
We denote $\tau_t=\max\{kh_N, kh_N\leq t\}.$
Multiplying \eqref{eq:Ito} by
$\frac1N\mathbf1_{\kappa_i^N>kh_N}\gamma(X_0^{N,i}),$
 summing over~$i$ and $k$, and integrating by parts, one obtains
\begin{align}\label{eq:apres_IPP}
\int_\R g_0H*\tilde\mu_0^N=&-\int_0^\infty\int_\R\partial_tg_tH*\tilde\mu_t^N\dx
t+\int_0^\infty\int_\R(-\Delta)^{\frac\alpha2}g_tH*\tilde\mu_t^N\dx t\nonumber\\
&+\frac1N\int_0^\infty\sum_{\kappa_i^N>\tau_t}\gamma(X_0^{N,i})
\partial_x\varphi_t(X_t^{N,i})A'\left(H*\tilde\mu_{kh_N}(X_{kh_N}^{N,i})\right)\dx t\nonumber\\
&+\frac1N\int_{(0,\infty)\times\R}\sum_{\kappa_i^N>\tau_t}\gamma(X_0^{N,i})\left(\varphi_t(X_{t-}^{N,i}+\sigma y)
-\varphi_t(X_{t-}^{N,i})\right)\tilde\nu^i(\dx y,\dx t)\nonumber\\
&-\frac1N\sum_{k=0}^\infty\sum_{\kappa_i^N=(k+1)h_N}\gamma(X_0^{N,i}))
\varphi_{(k+1)h_N}(X_{(k+1)h_N}^{N,i}),
\end{align}

Combining an adaptation of Lemma \ref{lem:transport}, stated in Section
\ref{sect:preuves_lemmes}, with $A$ replacing $\eta$, and integrating by
parts, the difference
$$\frac1N\int_0^\infty\sum_{\kappa_i^N>\tau_t}\gamma(X_0^{N,i})
\partial_x\varphi_t(X_t^{N,i})A'\left(H*\tilde\mu_{kh_N}(X_{kh_N}^{N,i})\right)\dx
t+\int_0^\infty\int_\R\partial_x g_tA(H*\tilde\mu_t^N)\dx
t$$vanishes in $\LL^1.$
Using an adaptation Lemma \ref{lem:martingale}, the
the fourth term in the right hand side of~\eqref{eq:apres_IPP}
vanishes in~$\LL^2.$ The fifth term vanishes in~$\LL^1$ since 
\begin{align*}
&\left|\frac1N\sum_{k=0}^\infty\sum_{\substack{\kappa_i^N=(k+1)h_N}}\gamma(X_0^{N,i})\varphi_{(k+1)h_N}(X_{(k+1)h_N}^{N,i})\right|\\
\leq&\frac1N\sum_{k=0}^\infty~~\sum_{\substack{\textrm{pairs }\{i,j\}\textrm{ killled}\\\textrm{at time }(k+1)h_N}}
\left|\varphi_{(k+1)h_N}\left(X_{(k+1)h_N}^{N,i}\right)-\varphi_{(k+1)h_N}\left(X_{(k+1)h_N}^{N,j}\right)\right|\\
\leq&K\varepsilon_N.
\end{align*}

As a consequence,
$\E^{\pi^N}|F(Q)|=\E|F(\mu^N)|$ vanishes. We conclude by regularizing
the function~$\gamma$ as in the two first points, that
$\E^{\pi^\infty}|F(Q)|=0$. Thus, $F(Q)=0$ almost surely, so that
$H*\tilde Q$ almost surely satisfies the weak formulation.

The case $\alpha=2$ is treated in the same way, the only difference
lying in the nature of the stochastic integral.

\end{enumerate}

\subsection{Proofs of Lemmas \ref{lem:martingale} to \ref{lem:controle_sauts}}\label{sect:preuves_lemmes}
In this section, we give the proofs of the previously admitted lemmas of
Section \ref{sect:preuve}.
\begin{proof}[Proof of Lemma \ref{lem:martingale}]
Since the particles are driven by
independent stable processes and since the inequality $|w_t^i|\leq\frac
KN$ holds for some constant $K$ not depending on $t$, $i$ and $N$,
\begin{align*}
\E M_N^2=&\E\left|\sum_{\kappa_i^N>0}\int_0^\infty w_{\tau_t}^i\chi_t^N\left(\int_{\{|y|\leq r_N\}}
\left(\varphi_t(X_{t-}^{N,i,r_N}+\sigma_Ny)-\varphi_t(X_{t-}^{N,i,r_N})\right)\tilde \nu^{i,r_N}(\dx y,\dx t)\right)\right|^2\\
\leq&2\sigma_N^2c_\alpha\E\left(\sum_{\kappa_i^N>0}\int_0^\infty
(w_{\tau_t}^i)^2\chi_t^N\int_{\{|y|\leq r_N\}}\left(y\|g_t\|_\infty\right)^2
\frac{\dx y\dx t}{|y|^{1+\alpha}}\right)\\
\leq&K\frac {\sigma_N^2r_N^{2-\alpha}}N\int_0^\infty\|g_t\|_\infty^2\dx t.
\end{align*}
A similar proof with stochastic integrals against Brownian motion yields the result for $\alpha=2.$
\end{proof}
\begin{proof}[Proof of Lemma \ref{lem:T2}]
Integrating by parts, one finds
$$\sum_{i=1}^N\int_0^\infty w_t^i\partial_t\varphi_t\left(X_t^{N,i,r_N}\right)\dx t
=-\int_0^\infty\int_\R\eta(H*\tilde\mu_t^{N,r_N})\partial_tg_t\dx t
+\int_0^\infty\int_\R\eta(\tilde\mu_t^{N,r_N}(\R))\partial_tg_t\dx t$$
yielding, from Lemma \ref{lem:transport} below,
$$\E\left|\sum_{i=1}^N\int_0^\infty w_{\tau_t}^i\partial_t\varphi_t\left(X_t^{N,i,r_N}\right)\dx  t
+\int_0^\infty\int_\R\eta(H*\tilde\mu_t^{N,r_N})\partial_tg_t\dx t
-\int_0^\infty\int_\R\eta(\tilde\mu_t^{N,r_N}(\R))\partial_tg_t\dx t\right|
\underset{N\rightarrow\infty}\longrightarrow0.$$
From the constancy of $\tilde\mu_t^{N,r_N}(\R)$ and an integration by parts, one has 
$$-T_N^1+\int_0^\infty\int_\R\eta(\tilde\mu_t^{N,r_N}(\R))\partial_tg_t\dx t=-\int_\R g_0\eta(H*\tilde\mu_0^{N,r_N}).$$
 Another integration by parts yields
\begin{align*}
&2c_\alpha\sum_{i=1}^N\int_0^\infty w_t^i\chi_t^N
\int_{\{|y|\leq r_N\}}\varphi_t\left(X_t^{N,i,r_N}+\sigma_Ny\right)
-\varphi_t\left(X_t^{N,i,r_N}\right)-\sigma_Ny\partial_x\varphi_t\left(X_t^{N,i,r_N}\right)\frac{\dx y\dx t}{|y|^{1+\alpha}}\\
=&-2c_\alpha\int_0^\infty\chi_t^N\int_{\{|y|\leq r_N\}}
\int_\R\left(g_t(x+\sigma_Ny)-g_t(x)-\sigma_Ny\partial_xg_t(x)\right)\eta(H*\tilde\mu_t^{N,r_N}(x))\frac{\dx
  x\dx y\dx t}{|y|^{1+\alpha}}\\
&+2c_\alpha\int_0^\infty\chi_t^N\eta(\tilde\mu_t^{N,r_N}(\R))\int_{\{|y|\leq r_N\}}
\int_\R\left(g_t(x+\sigma_Ny)-g_t(x)-\sigma_Ny\partial_xg_t(x)\right)\frac{\dx
  x\dx y\dx t}{|y|^{1+\alpha}}\\
=&-2c_\alpha\int_0^\infty\chi_t^N\int_{\{|y|\leq r_N\}}
\int_\R\left(g_t(x+\sigma_Ny)-g_t(x)-\sigma_Ny\partial_xg_t(x)\right)\eta(H*\tilde\mu_t^{N,r_N}(x))\frac{\dx
  x\dx y\dx t}{|y|^{1+\alpha}}.
\end{align*}
Moreover, from the regularity of $A$ and~$\eta$, it holds
$$w_{\tau_t}^iA'\left(H*\tilde\mu_{\tau_t}^{N,r_N}(X_{\tau_t}^{N,i,r_N})\right)=
\psi\left(H*\tilde\mu_{\tau_t}^{N,r_N}(X_{\tau_t}^{N,i,r_N})\right)
-\psi\left(H*\tilde\mu_{\tau_t}^{N,r_N}(X_{\tau_t}^{N,i,r_N}-)\right)+o\left(\frac1N\right),$$
so that
$$\E\left|2\sum_{i=1}^N\int_0^\infty
w_{\tau_t}^i\chi_t^N\partial_x\varphi_t\left(X_t^{N,i,r_N}\right)
A'\left(H*\tilde\mu_{\tau_t}^{N,r_N}(X_{\tau_t}^{N,i,r_N})\right)\dx
t+2\int_0^\infty\chi_t^N\int_\R\partial_xg_t\psi(H*\tilde\mu_t^{N,r_N})\dx
t\right|\underset{N\to\infty}\to0,$$
from an adaptation of Lemma \ref{lem:transport} (replacing $\eta$ by
$\psi$ in the definition of $w_t^i$). This concludes the proof of the
first item of Lemma \ref{lem:T2}.

To prove the second item, observe that the change of variable
$z=\sigma_Ny$ yields, for~$r_N\leq\frac1{\sigma_N}$,
\begin{align*}
&\left|2c_\alpha\int_0^\infty\chi_t^N\int_\R\int_{\{|y|\leq
  r_N\}}\eta(H*\tilde\mu_t^{N,r_N}(x))\left(g_t(x+\sigma_Ny)-g_t(x)-\sigma_Ny\partial_xg_t(x)\right)\frac{\dx y\dx
  x\dx t}{|y|^{1+\alpha}}\right|\\\leq&2c_\alpha\sigma_N^\alpha\int_0^\infty\chi_t^N\int_\R\int_{\{|z|\leq
  1\}}\left|\eta(H*\tilde\mu_t^{N,r_N}(x))\left(g_t(x+z)-g_t(x)-z\partial_xg_t(x)\right)\right|\frac{\dx z\dx
  x\dx t}{|z|^{1+\alpha}}.
\end{align*}
\end{proof}

\begin{proof}[Proof of Lemma \ref{lem:T3}]
First notice that
$$T_N^3=\sum_{\kappa_i^N>0}\int_0^\infty(1-\chi_t^N)\int_{\{|y|>r_N\}}\left(\int_\R\varphi_t\dx\rho_{t-}^{y,i}\right)\nu^{i,r_N}(\dx
y,\dx t),$$
with $\rho$ defined by the following formula: ($\bar\mu_t^{y,i,N,r_N}$ being the measure obtained by
moving in the expression of $\tilde\mu_t^{N,r_N}$ the particle
$X^{N,i,r_N}_t$ to the position
$X^{N,i,r_N}_t+\sigma_Ny$)
$$\rho_t^{y,i}=\partial_x\left(\eta(H*\bar\mu^{y,i,N,r_N}_t)-\eta(H*\tilde\mu^{N,r_N}_t)\right).$$

To prove the second item in Lemma \ref{lem:T3}, we integrate by parts,
and, using the definition 
of~$\bar\mu^{y,i,N,r_N}$ and the compactness of the support of~$g$, it holds
\begin{equation}\label{eq:majo}
\left|\int_\R\varphi_t\dx\rho_t^{y,i}\right|
=\left|\int_\R g_t\left(\eta(H*\bar\mu_t^{y,i,N,r_N})-\eta(H*\tilde\mu_t^{N,r_N})\right)\right|\\
\leq K\frac{(\sigma_Ny)\wedge 1}N,
\end{equation}
so that 
$$\E|T_N^3|\leq K\int_0^\infty(1-\chi_t^N)
\int_{\{|y|>r_N\}}(\sigma_N y)\wedge 1\frac{\dx
y\dx t}{|y|^{1+\alpha}}\leq K(\sigma_N^\alpha+\sigma_Nr_N^{1-\alpha}).$$

Now let us prove the first item of Lemma \ref{lem:T3}.
Applying the same martingale argument as the one used to
 prove $\E|M_N|\to0,$ and using the upper bound $K/N$ in \eqref{eq:majo}, one has
$$\E\left|T_N^3-2c_\alpha\int_0^\infty(1-\chi_t^N)\int_{\{|y|>r_N\}}\left(\sum_{\kappa_i^N>t}\int_\R\varphi_t\dx
\rho_t^{y,i}\right)\frac{\dx y\dx
t}{|y|^{1+\alpha}}\right|^2\leq \frac K{r_N^\alpha N}.$$

Let us give a more explicit expression for $\rho_t^{y,i}.$ For simplicity, we
denote
$$\tilde w_t^i=\mathbf1_{\kappa_i^N>t}\left[\eta\left(\frac1N\sum_{\underset{\kappa_j^N>t}{j\neq
    i}}\gamma(X_0^j)\mathbf1_{X^{N,j,r_N}_t\leq
    X^{N,i,r_N}_t+\sigma_Ny}+\frac{\gamma(X_0^i)}N\right)-
\eta\left(\frac1N\sum_{\underset{\kappa_j^N>t}{j\neq i}}\gamma(X_0^j)\mathbf1_{X^{N,j,r_N}_t\leq
    X^{N,i,r_N}_t+\sigma_Ny}\right)\right]$$
and for $i\neq j,$
$$\tilde w_t^{i,j,\pm}=\mathbf1_{\kappa_i^N>t}\mathbf1_{\kappa_j^N>t}
\left[\eta\left(\frac1N\sum_{\underset{\kappa_k^N>t}{k\neq j}}\gamma(X_0^k)\mathbf1_{X^{N,k,r_N}_t\leq
    X^{N,i,r_N}_t}\pm\frac{\gamma(X_0^j)}N\right)-\eta\left(\frac1N\sum_{\underset{\kappa_k^N>t}{k\neq j}}\gamma(X_0^k)\mathbf1_{X^{N,k,r_N}_t\leq
    X^{N,i,r_N}_t}\right)\right].$$
One can write
\begin{align}\label{eq:rho}
\rho_t^y:=\sum_{\kappa_i^N>t}\rho_t^{y,i}=&\sum_{\kappa_i^N>t} \tilde
w_t^i\delta_{X^{N,i,r_N}_t+\sigma_Ny}-
\sum_{\kappa_i^N>t}w_t^i\delta_{X^{N,i,r_N}_t}\nonumber\\
&+\sum_{\kappa_i^N>t}\left(\sum_{\kappa_j^N>t}\left(\tilde
    w_t^{i,j,+}-w_t^i\right)\mathbf1_{X_t^{N,i,r_N}<X_t^{N,j,r_N}}
\mathbf1_{X^{N,j,r_N}_t+\sigma_Ny<X^{N,i,r_N}_t}\right)\delta_{X^{N,i,r_N}_t}\nonumber\\
&+\sum_{\kappa_i^N>t}\left(\sum_{\kappa_j^N>t}\left(\tilde w_t^{i,j,-}-w_t^i\right)
\mathbf1_{X_t^{N,j,r_N}<X_t^{N,i,r_N}}\mathbf1_{X^{N,i,r_N}_t<X^{N,j,r_N}_t+\sigma_Ny}\right)\delta_{X^{N,i,r_N}_t}.
\end{align}
In this expression, the two first terms deal with particles
jumping from the site $X_t^{N,i,r_N}$ to the site
$X_t^{N,i,r_N}+\sigma_Ny$,
while the third term corresponds to the jump from right to left of the
particle labelled~$j$ above the particle labelled $i$ and, conversely,
the fourth term corresponds to the jumps of particle $j$ from left to
right over particle $i$.
Notice that this last equality, as well as \eqref{eq:sigma} below, only holds when each $X_t^{N,i,r_N}+\sigma_Ny$ is
distinct from all~$X_t^{N,j,r_N}.$ However, for all $t$, this condition holds $\dx
y$-almost everywhere, which is enough for our purpose.

In the entropic formulation \eqref{eq:entropifr}, the term that should
appear for large jumps is given by 
$$2c_\alpha\int_0^\infty\int_{\{|y|>r_N\}}\left(\int_\R\varphi_t\dx \sigma_t^y\right)\frac{\dx y\dx t}{|y|^{1+\alpha}},$$
where
\begin{align}\label{eq:sigma}
\sigma_t^y=&\partial_x\left(\eta'(H*\tilde\mu^{N,r_N}_t)
\left(H*\tilde\mu^{N,r_N}_t(\cdot-\sigma_Ny)-H*\tilde\mu^{N,r_N}_t\right)\right)\nonumber\\
=&\frac1N\sum_{\kappa_i^N>t}\gamma(X_0^i)\eta'\left(H*\tilde\mu^{N,r_N}_t(X^{N,i,r_N}_t+\sigma_Ny)\right)\delta_{X^{N,i,r_N}_t+\sigma_Ny}
-\frac1N\sum_{\kappa_i^N>t}\gamma(X_0^i)\eta'\left(H*\tilde\mu_t^{N,r_N}(X_t^{N,i,r_N})\right)\delta_{X^{N,i,r_N}_t}\nonumber\\
&+\sum_{\kappa_i^N>t}\left(H*\tilde\mu^{N,r_N}_t(X^{N,i,r_N}_t-\sigma_Ny)-H*\tilde\mu_t^{N,r_N}(X_t^{N,i,r_N})\right)\nonumber\\
&~~~~~~~~~~~\times\left(\eta'\left(H*\tilde\mu_t^{N,r_N}(X_t^{N,i,r_N})\right)
-\eta'\left(H*\tilde\mu_t^{N,r_N}(X_t^{N,i,r_N}-)\right)\right)\delta_{X_t^{N,i,r_N}}.
\end{align}
When computing the difference $\rho_t^y-\sigma_t^y$ integrated against
some bounded function, using Taylor expansions for $\eta,$ one can
check that, up to an error term of order $\mathcal O\left(\frac1N\right)$ 
the first terms in the right hand side of \eqref{eq:rho} and \eqref{eq:sigma}
cancel each other, the second terms as well, and so does the sum of the
two last term in \eqref{eq:rho} with the
last one in~\eqref{eq:sigma}.
Consequently, 
$$\left|\int_0^\infty(1-\chi_t^N)\int_{\{|y|>r_N\}}
\left(\int_\R\varphi_t\dx\rho_t^y\right)\frac{\dx y\dx t}{|y|^{1+\alpha}}
-\int_0^\infty(1-\chi_t^N)\int_{\{|y|>r_N\}}\left(\int_\R\varphi_t\dx
  \sigma_t^y\right)
\frac{\dx y\dx t}{|y|^{1+\alpha}}\right|\leq \frac K{Nr_N^\alpha}.$$
This concludes the proof.
\end{proof}

\begin{proof}[Proof of Lemma \ref{lem:T4}]
For a time $\zeta_m$ of the form $kh_N,$ no particle moved in the
interval $(\zeta_{m-1},\zeta_m)$, so that
$w_{\zeta_m}^i-w_{\zeta_{m-1}}^i=0,$ unless the particle labelled $i$ has
been killed at time $\zeta_m.$ Hence,
\begin{align*}T_N^4=&\sum_{i=1}^N\sum_{\substack{\zeta_m\textrm{ of the}\\
\textrm{form }kh_N}}(w_{\zeta_m}^i-w_{\zeta_{m-1}}^i)\varphi_{\zeta_m}(X_{\zeta_m}^{N,i,r_N})\\
=&-\sum_{\substack{\zeta_m\textrm{ of the}\\
\textrm{form }kh_N}}~~\sum_{\kappa_i^N=\zeta_m}w_{\zeta_{m-1}}^i\varphi_{\zeta_m}(X_{\zeta_m}^{N,i,r_N}).
\end{align*}
This sum is actually a sum over pairs of close particles with opposite
signs, thus
\begin{align*}
|T_N^4|&=\left|\sum_{\substack{\zeta_m\textrm{ of the}\\\textrm{form }kh_N}}~~
\sum_{\substack{\textrm{pairs }\{i,j\}\textrm{ of particles}\\\textrm{killed at time }\zeta_m}}
\left(w_{\zeta_{m-1}}^i\varphi_{\zeta_m}(X_{\zeta_m}^{N,i,r_N})
+w_{\zeta_{m-1}}^j\varphi_{\zeta_m}(X_{\zeta_m}^{N,j,r_N})\right)\right|\\
&\leq\sum_{\substack{\zeta_m\textrm{ of the}\\\textrm{form }kh_N}}~~
\sum_{\substack{\textrm{pairs }\{i,j\}\textrm{ of particles}\\\textrm{killed at time }\zeta_m}}
\left|w_{\zeta_{m-1}}^i+w_{\zeta_{m-1}}^j\right|\|\varphi\|_\infty
+\left|w_{\zeta_{m-1}}^j\right|\left|\varphi_{\zeta_m}(X_{\zeta_m}^{N,i,r_N})-\varphi_{\zeta_m}(X_{\zeta_m}^{N,j,r_N})\right|\\
&\leq
 K\left(\frac1N+\varepsilon_N\right).
\end{align*}
Indeed, a couple $(i,j)$ of killed particles is such that
$|X_{\zeta_m}^{N,i,r_N}-X_{\zeta_m}^{N,j,r_N}|\leq\varepsilon_N$ and is
made of particles with opposite signs, so that
$$|w_{\zeta_{m-1}}^i+w_{\zeta_{m-1}}^j|=\left|\left(\gamma(X_0^i)+\gamma(X_0^j)\right)\eta'(H*\tilde\mu_t^{N,r_N}(X_t^{N,i,r_N}))+\mathcal
O\left(\frac1{N^2}\right)\right|\leq \frac K{N^2}.$$

\end{proof}

\begin{proof}[Proof of Lemma \ref{lem:controle_sauts}]
Notice that from independence of
the increments,
denoting by $\mathcal L^{\leq r}$ a L\'evy process with L\'evy measure
$c_\alpha\mathbf1_{|y|\leq r}\frac{\dx y}{|y|^{1+\alpha }}$, it holds
\begin{align*}
\PP(B_N)&=\PP(\sigma_N|\mathcal L_{h_N}^{\leq r_N}|\leq\varepsilon_N/4)^{NT/h_N}\\
&=\left(1-\PP\left(\sigma_Nr_N|\mathcal L^{\leq1}_{h_Nr_N^{-\alpha}}|\geq\varepsilon_N/4\right)\right)^{NT/h_N}.
\end{align*}
Since the L\'evy measure $c_\alpha\mathbf1_{|y|\leq1}\frac{\dx y}{|y|^{1+\alpha}}$ has compact support, the
random variables $\mathcal L_t^{\leq1}$ have exponential moments, and Chernov's
inequality yields
$$\PP\left(\sigma_Nr_N|\mathcal L^{\leq1}_{h_Nr_N^{-\alpha}}|\geq\varepsilon_N/4\right)\leq
\mathbb E\left(e^{\left|\mathcal L^{\leq1}_{h_Nr_N^{-\alpha}}\right|}\right)e^{-\varepsilon_N/4\sigma_Nr_N}
=e^{Kh_Nr_N^{-\alpha}-\varepsilon_N/4\sigma_Nr_N},$$where the constant $K$ does
not depend on $N$.

In the Brownian case $\alpha=2$, we use the tail estimate
$\int_M^\infty e^{-x^2}\dx x\leq Ke^{-M^2}$ for positive $M$.
\end{proof}

\begin{lem}\label{lem:max_somme}
Let $a_1\leq\hdots\leq a_N$ and $b_1\leq\hdots\leq b_N$ be two nondecreasing sequences of reals
numbers. Then the quantity $\sum_{i=1}^Na_ib_{\sigma(i)}$ for some
permutation $\sigma$ is maximal when $\sigma(i)=i$ for all~$i$.
\end{lem}
\begin{proof}
From optimal transportation theory (see \cite[page 75]{villani-03}), the quantity
$\sum_{i=1}^N(a_i-b_{\sigma(i)})^2$ is minimal when~$\sigma$ is the
identity. Expanding the square, we see that
$\sum_{i=1}^N(a_i-b_{\sigma(i)})^2=\sum_{i=1}^N(a_i^2+b_i^2)-2\sum_{i=1}^Na_ib_{\sigma(i)}.$
Thus, $\sum_{i=1}^Na_ib_{\sigma(i)}$ is maximal if and only if
$\sum_{i=1}^N(a_i-b_{\sigma(i)})^2$ is minimal, concluding the proof.
\end{proof}

\begin{lem}\label{lem:transport}
Let $f$ be some bounded function with
compact support on $[0,\infty)\times\R$ which is smooth with respect to
the space variable. If $h_N$ vanishes and $\sigma_N$ is bounded, it holds
$$\lim_{N\to\infty}\E\left|\sum_{\kappa_i^N>0}\int_0^\infty\left(w_t^i-w_{\tau_t}^i\right)
f_t\left(X_t^{N,i,r_N}\right)\dx t\right|=0.$$
\end{lem}
\begin{proof}
First notice that when $t$ is not in an interval $[kh_N,(k+1/2)h_N]$, it holds
$w_t^i=w_{\tau_t}^i,$ since no particle moved between $\tau_t$ and $t.$
Then, one can write, from the assumptions on $f$,
\begin{align*}
\left|\sum_{\kappa_i^N>0}\int_0^\infty\left(w_t^i-w_{\tau_t}^i\right)
f_t\left(X_t^{N,i,r_N}\right)\dx
t\right|
\leq&\left|\sum_{\kappa_i^N>0}\int_0^T \chi_t^N
\left(  w_t^if_t(X_t^{N,i,r_N})-w_{\tau_t}^if_t(X_{\tau_t}^{N,i,r_N})\right)\dx
  t\right|\\
&
+\frac
KN\sum_{\kappa_i^N>0}\int_0^T\chi_t^N\left|X_t^{N,i,r_N}-X_{\tau_t}^{N,i,r_N}\right|\wedge1\dx
t.
\end{align*}
Integrating by parts, it holds:
\begin{align*}
&\left|\sum_{\kappa_i^N>0}\int_0^T\chi_t^N\left(w_t^if_t(X_t^{N,i,r_N})-w_{\tau_t}^i
f_t(X_{\tau_t}^{N,i,r_N})\right)\dx t\right|\\
=&\left|\int_0^T\chi_t^N\int_\R\left(\eta\left(H*\tilde\mu_t^{N,r_N}(x)\right)
-\eta\left(H*\tilde\mu_{\tau_t}^{N,r_N}(x)\right)\right)\partial_xf_t(x)\dx
x\dx  t\right|\\
\leq&\frac
KN\int_0^T\chi_t^N\int_\R\left(\sum_{\kappa_i^N>0}\mathbf1_{X_t^{N,i,r_N}\leq
  x<X_{\tau_t}^{N,i,r_N}}+\mathbf1_{X_{\tau_t}^{N,i,r_N}\leq
  x<X_t^{N,i,r_N}}\right)\partial_xf_t(x)\dx x\dx t\\
\leq&\frac
KN\int_0^T\chi_t^N\sum_{\kappa_i^N>0}\left|X_t^{N,i,r_N}-X_{\tau_t}^{N,i,r_N}\right|\wedge1\dx t
\end{align*}
We conclude the proof by writing
\begin{align*}
\E\frac 1N\int_0^T\chi_t^N\sum_{\kappa_i^N>0}
\left|X_t^{N,i,r_N}-X_{\tau_t}^{N,i,r_N}\right|\wedge1\dx t&=
\E\int_0^T\chi_t^N\mathbf1_{\kappa_1^N>0}\left|X_t^{N,1,r_N}-X_{\tau_t}^{N,1,r_N}\right|\wedge1\dx
t\\
&\leq T\left(h_N\sup_{[-1,1]}|A'|+\E\left((\sigma_N|\Lambda^{N,1,r_N}_{h_N}|)\wedge1\right)\right).
\end{align*}
This last quantity vanishes when $h_N$ goes to 0.
\end{proof}

\section{Numerical results}

In this section, we illustrate our convergence results by some numerical
simulations. We simulated the solution to the fractional and the inviscid Burgers equations
$$\partial_tu+\frac12\partial_x(u^2)+\sigma^\alpha(-\Delta)^{\frac\alpha2}=0~\textrm{ and }~\partial_tu+\frac12\partial_x(u^2)=0,$$
corresponding to the choice $A(x)=x^2/2$, 
with different values for the parameter $\alpha.$

One can find an explicit exact
solution to the inviscid Burgers equation (see \cite{lax-57}) and we compare the result of
the simulation to this exact solution in the vanishing viscosity
setting. However, to our knowledge, no explicit
solutions exist in the case of a positive viscosity coefficient for $\alpha<2$, so that we have to compare 
the result of our
simulation with the one given by another numerical method. Here, we use a deterministic method,
introduced by Droniou in \cite{droniou-09}.

\subsection{Constant viscosity ($\sigma_N=\sigma$)}

\begin{figure}[!h]
\centerline{\epsfig{file=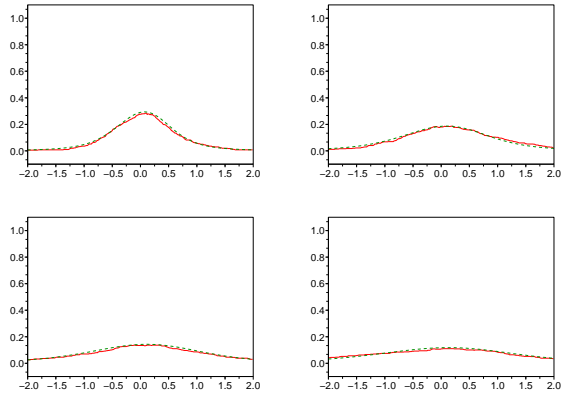,width=80mm,angle=0}}
\caption{\small Approximation of the conservation law with index $\alpha=1.5$.}\label{fig:alpha_grand_visq}
\end{figure}

\begin{figure}[!h]
\centerline{\epsfig{file=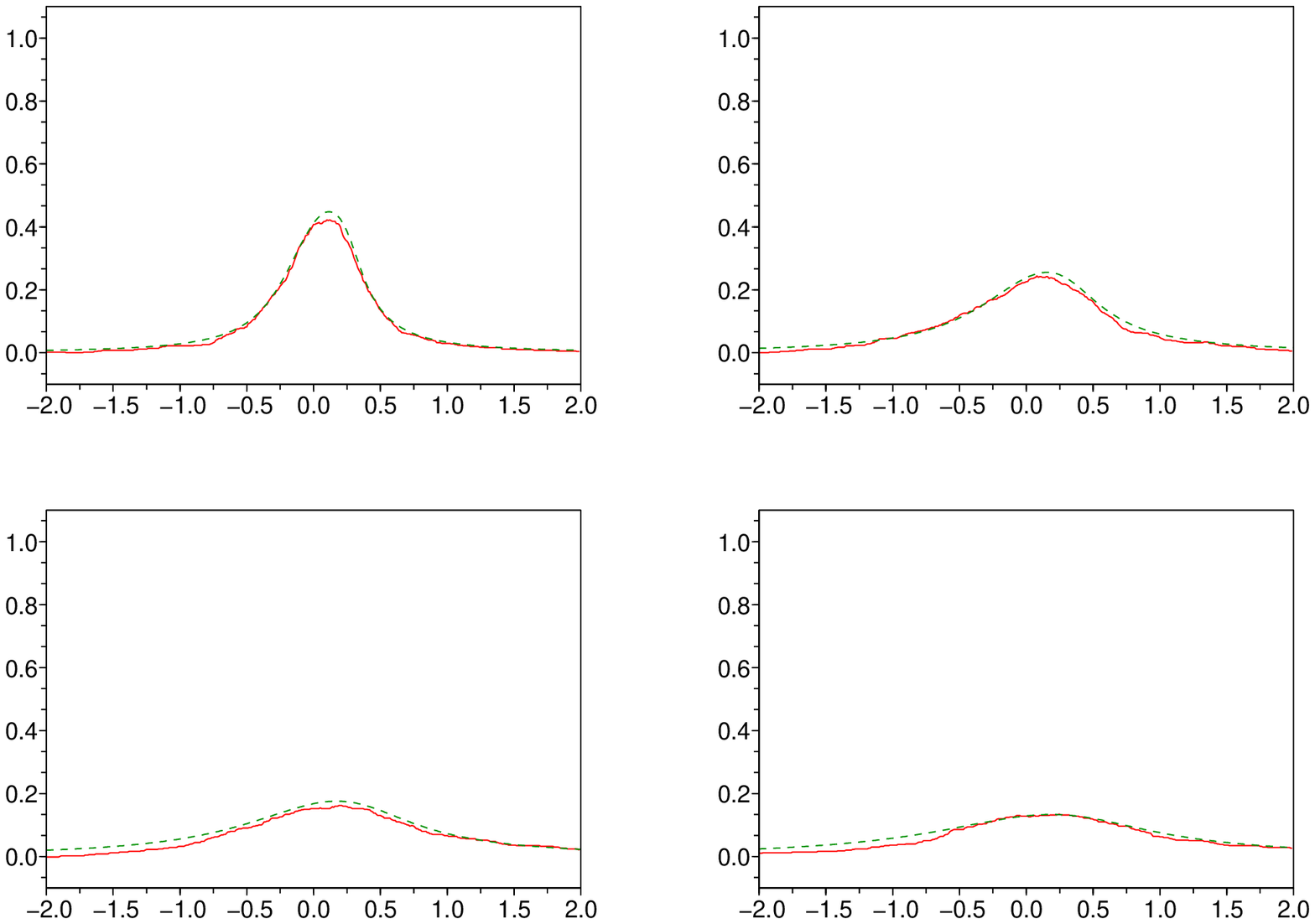,width=80mm,angle=0}}
\caption{\small Approximation of the conservation law with index $\alpha=1$.}\label{fig:alpha_moyen_visq}
\end{figure}

\begin{figure}[!h]
\centerline{\epsfig{file=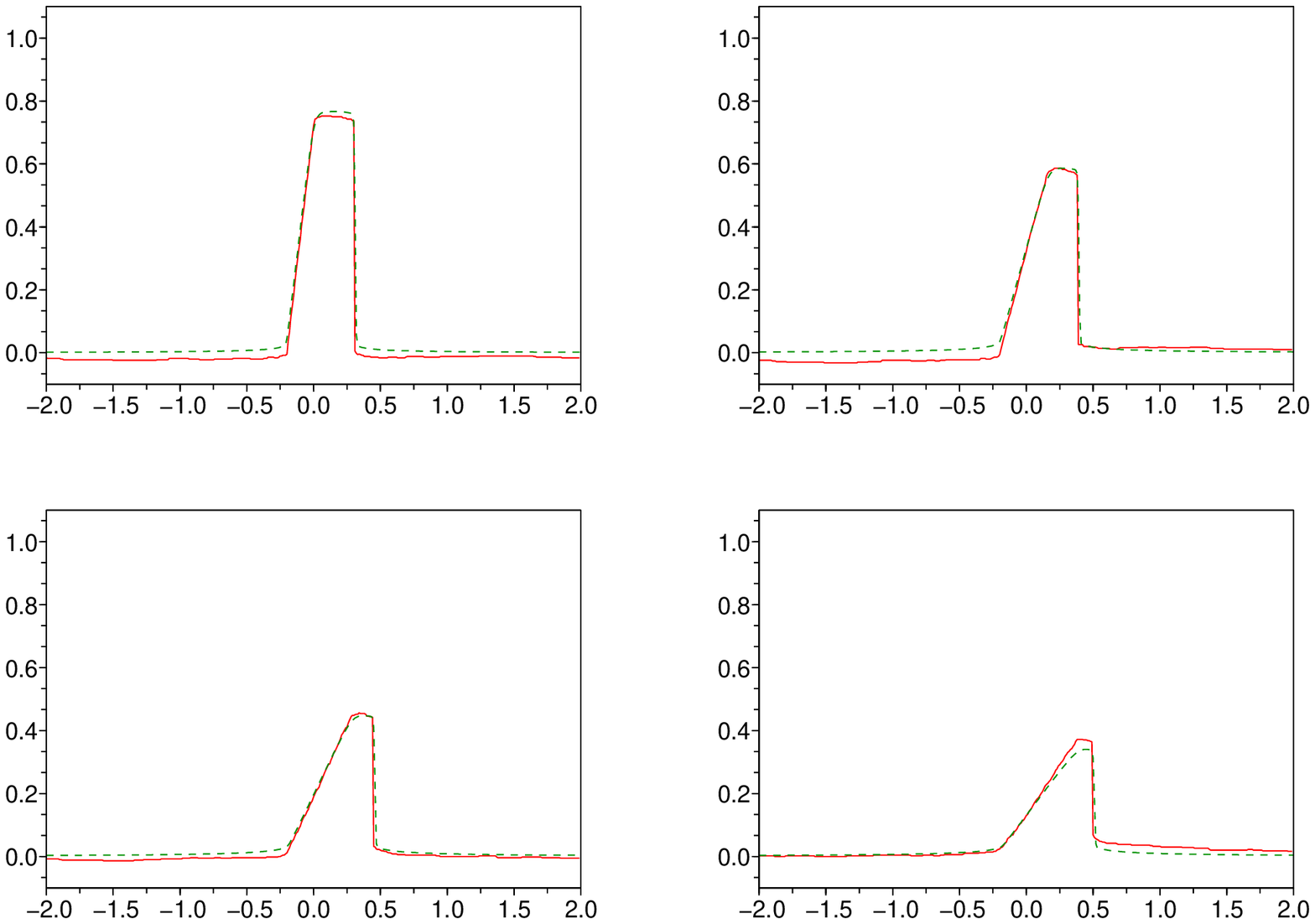,width=80mm,angle=0}}
\caption{\small Approximation of the conservation law with index $\alpha=0.1$.}\label{fig:alpha_petit_visq}
\end{figure}

We give three examples of approximation to the viscous conservation law.
On Figures \ref{fig:alpha_grand_visq}, \ref{fig:alpha_moyen_visq} and \ref{fig:alpha_petit_visq}, we show the approximation of the
viscous conservation law with respective index $\alpha=1.5$, $\alpha=1$ and $\alpha=0.1$ and diffusion
coefficient $\sigma=1$ using $N=1000$ particles, with parameters
$h=0.01$ and $\varepsilon=0.04$ at simulation times 0.25, 0.5, 0.75 and 1. The continuous line is the simulated
solution, and the dotted line is the ``exact'' solution obtained with the
determistic scheme of \cite{droniou-09} using small time and space steps.

\begin{figure}[!h]
\centerline{\epsfig{file=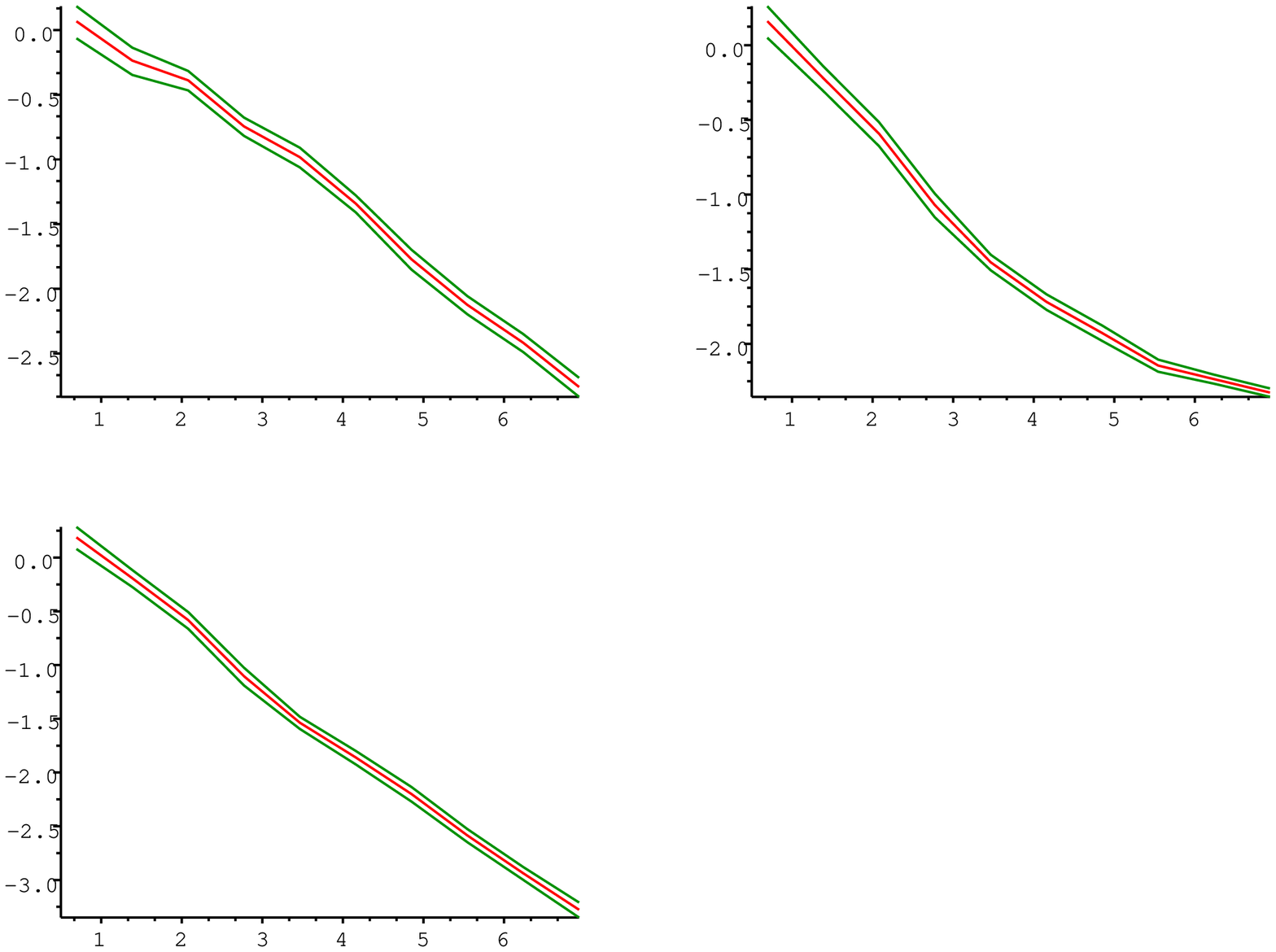,width=80mm,angle=0}}
\caption{\small Logarithmic error in the approximation of the conservation law
  with index $\alpha=0.5$, $1$ and $1.5$. The respective slopes are $-0.46$, $-0.41$ and $-0.56$.}
\label{fig:erreur_N}
\end{figure}

We now investigate the vanishing rate of the error, that is the Riemann
sum on the discretization grid associated to the integral in Theorems
\ref{theo:conv_alpha<1}, \ref{theo:conv_sigma0} and
\ref{theo:conv_alpha>1}.
On Figure \ref{fig:erreur_N} is depicted the logarithmic plot
of the error as a function of $N$ where we used the relation $h_N=10/N$,
and $\varepsilon_N=40/N$, with~$N$ ranging from 10 to 10000, in the three cases $\alpha=0.5$, $1$ and $1.5$.
In the case $\alpha<1$, this relation between $N$, $h_N$ and $\varepsilon_N$ satisfies the condition 
of Theorem~\ref{theo:conv_alpha<1}.
These pictures make us expect a convergence rate of $\frac1{\sqrt N}$,
corresponding to the optimal rate analyzed theoretically in~\cite{bossy-talay-96,bossy-talay-97}, in the
case $\alpha=2,$ without killing.

\subsubsection{Behaviour as $h\to0$}

\begin{figure}[!h]
\centerline{\epsfig{file=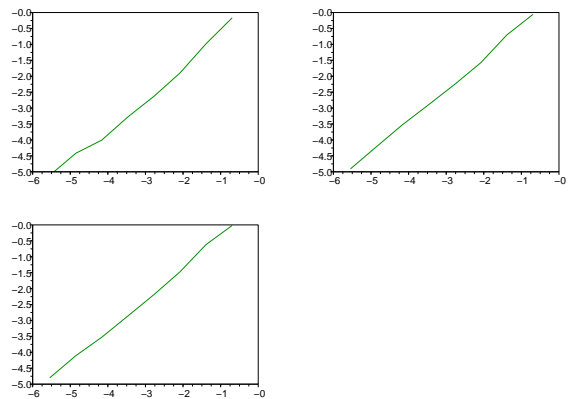,width=80mm,angle=0}}
\caption{\small Logarithmic plot of the error as $h$ vanishes, with a fixed
  number of particles, at respectively $\alpha=0.5$, $1$ and $1.5$. The slopes are
  equal to 1 up to an error of 0.01.}\label{fig:erreur_h}
\end{figure}

We give in Figure \ref{fig:erreur_h} the approximation error at fixed
number of particle, with a vanishing time step~$h$, in logarithmic plot. We set the parameter
$\varepsilon$ to be equal to $4h$ so that the condition of Theorem~\ref{theo:conv_alpha<1} is satified. We took $N=340000$ and
$\sigma=1$. We set $\alpha=0.5$, $\alpha=1$ and $\alpha=1.5$ respectively. The different parameters $h$ range
from $1$ to $2^{-8}$.
In \cite{bossy-talay-96,bossy-talay-97} it is shown, in case
$\alpha=2$ and the initial condition is monotonic, that the error is of order
$h.$ In view of Figure \ref{fig:erreur_h}, it seems that the convergence
rate is still of order $h$, even for $\alpha<2$ and any initial
condition with bounded variation.

\subsection{Vanishing viscosity ($\sigma_N\to0$)}

We consider the Burgers equation $$\partial_t v=\partial_x(u^2/2)$$ with
initial condition
$u_0(x)=\mathbf1_{[-3,-2]}-\mathbf1_{[2,3]},$
 which is the cumulative distribution function of the measure
 $\delta_{-3}-\delta_{-2}+\delta_2-\delta_3$. In that case, the
 solution of the Burgers equation is explicit and given by the
 expression
$$u(t,x)=\min\left(\frac{x+3}t,1\right)\mathbf1_{[-3,\min(-2+\frac
  t2,-3+\sqrt{2t},0)]}+\max\left(\frac{x-3}t,-1\right)\mathbf1_{[\max(2-\frac t2,3-\sqrt{2t},0),3]}.$$
We compare the function $u$ to the function obtained by running the
Euler scheme with a small diffusion coefficient $\sigma.$ One can expect the approximation to be
better for large values of $\alpha.$ Indeed, for small values of $\alpha$,
the particles tend to jump very far away, and subsequently
``disappear'' from the simulation. The consequence of this behavious is
that the solution is somehow decreased by a multiplicative constant.

\begin{figure}[!h]
\centerline{\epsfig{file=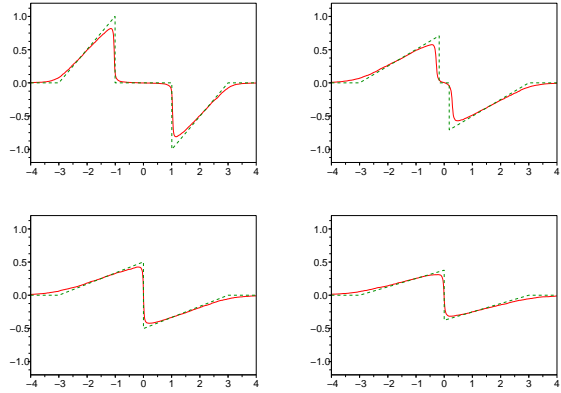,width=80mm,angle=0}}
\caption{\small Approximation of the inviscid conservation law by a fractional
Euler scheme with index $\alpha=1.5$ and diffusion coefficient $0.1$.}\label{fig:alpha_grand}
\end{figure}

\begin{figure}[!h]
\centerline{\epsfig{file=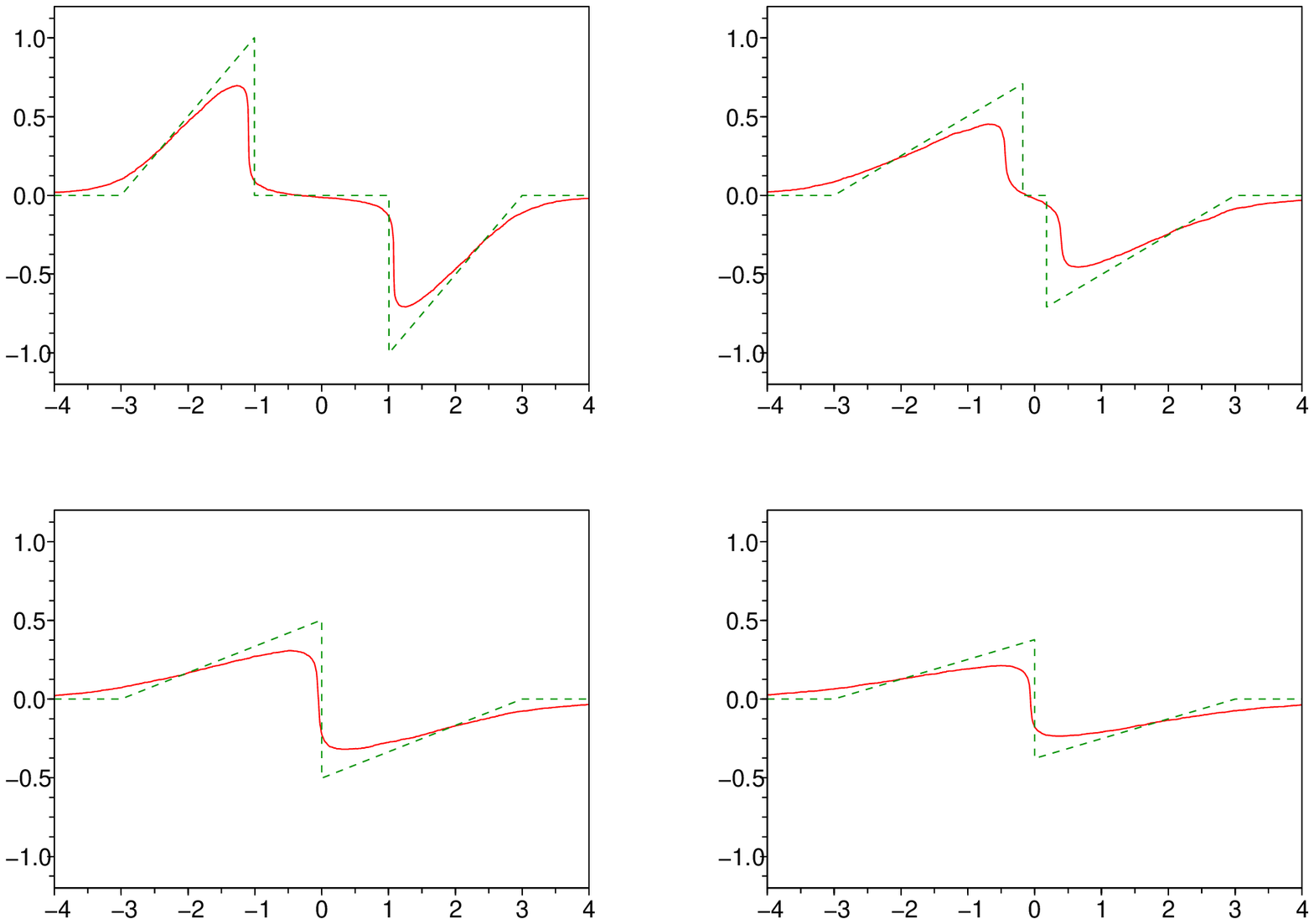,width=80mm,angle=0}}
\caption{\small Approximation of the inviscid conservation law by a fractional
Euler scheme with index $\alpha=1$ and diffusion coefficient
$0.1$.}\label{fig:alpha_moyen}
\end{figure}

\begin{figure}[!h]
\centerline{\epsfig{file=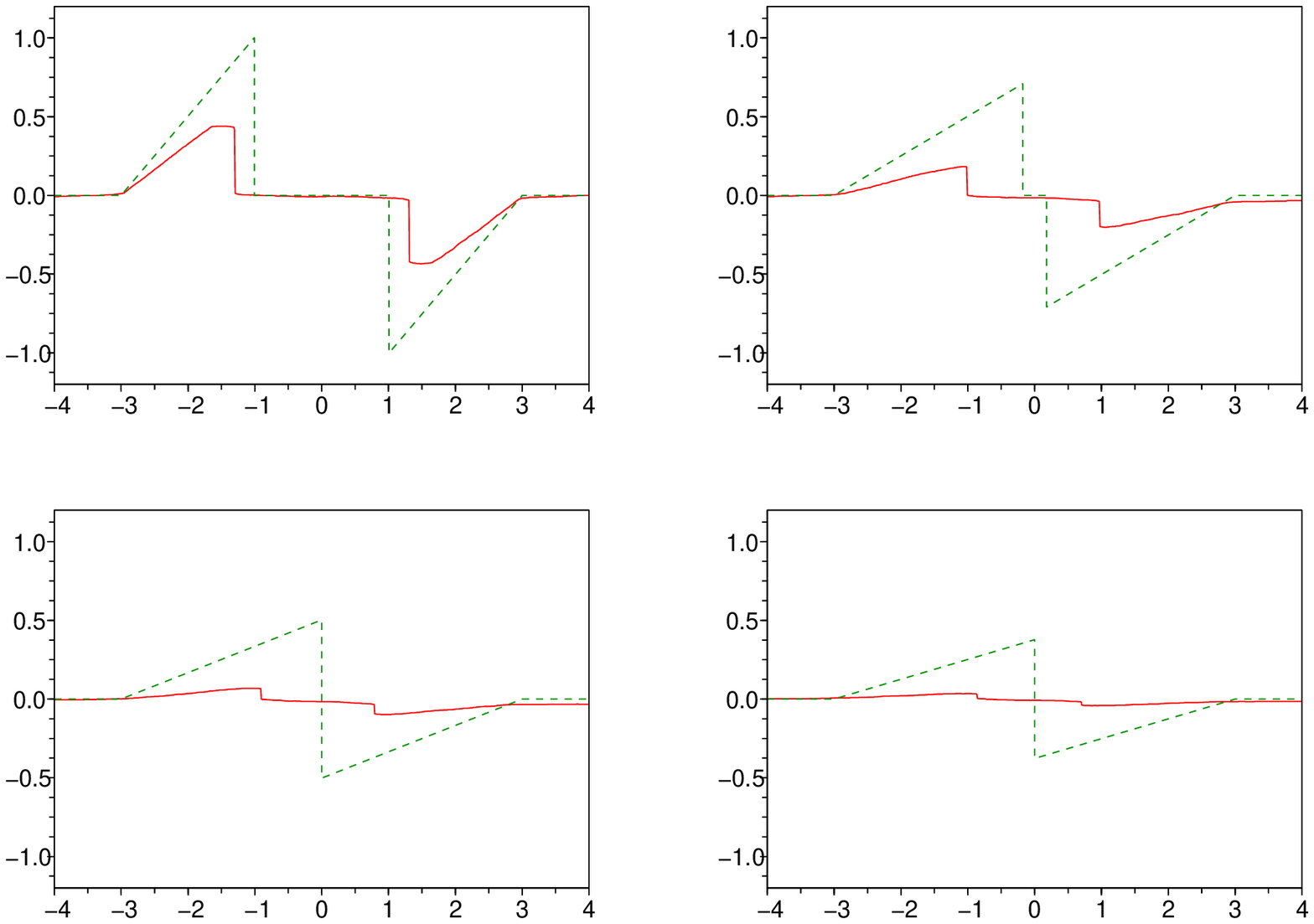,width=80mm,angle=0}}
\caption{\small Approximation of the inviscid conservation law by a fractional
Euler scheme with index $\alpha=0.1$ and diffusion coefficient $10^{-4}$.}\label{fig:alpha_petit_mauvais}
\end{figure}

\begin{figure}[!h]
\centerline{\epsfig{file=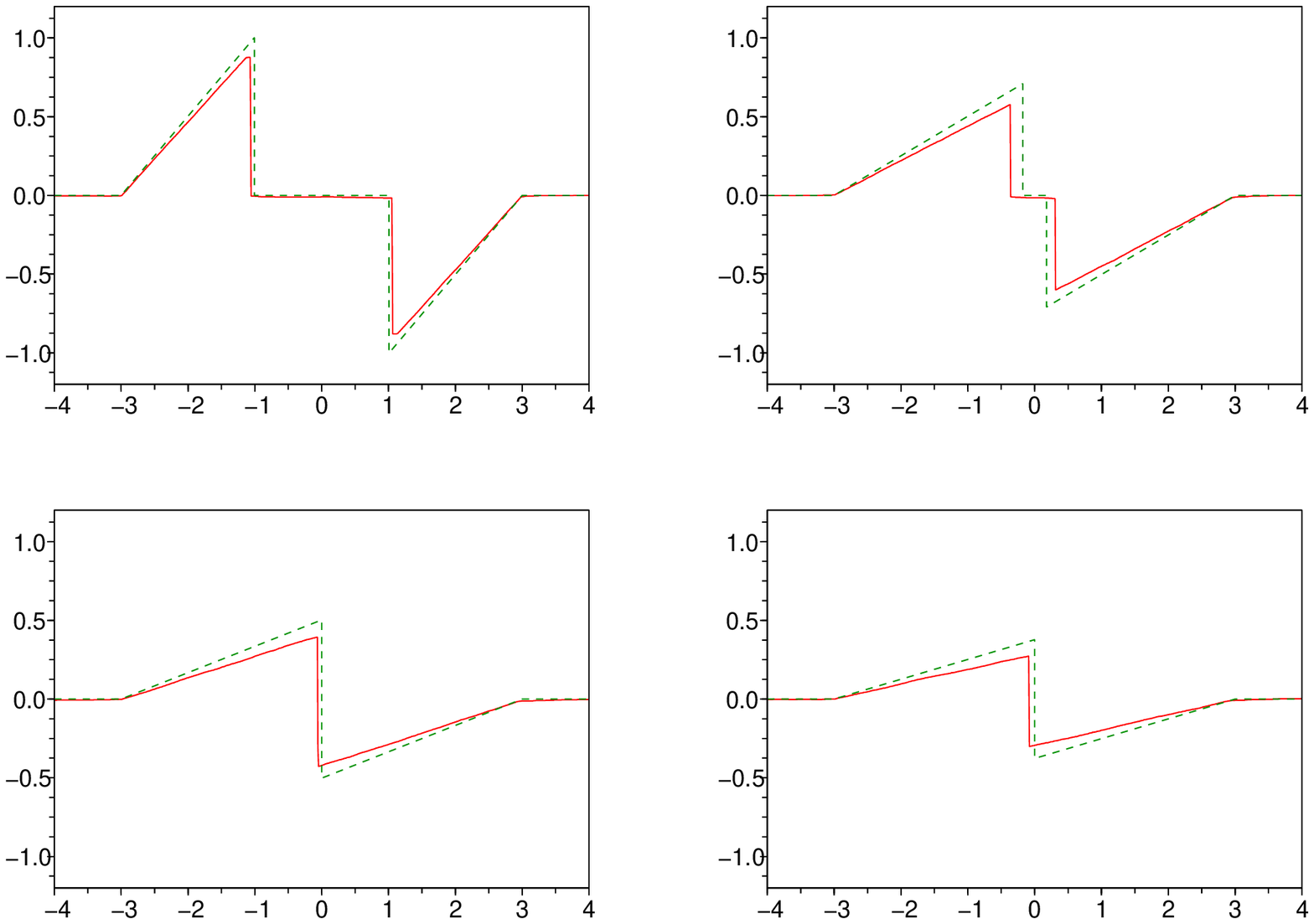,width=80mm,angle=0}}
\caption{\small Approximation of the inviscid conservation law by a fractional
Euler scheme with index $\alpha=0.1$ and diffusion coefficient $10^{-12}$.}\label{fig:alpha_petit_bon}
\end{figure}

For large values of $\alpha$, the approximation is quite good, even for not
so small diffusion coefficients. Figure \ref{fig:alpha_grand} gives the
result of the simulation of the Euler scheme with parameters
$\alpha=1.5$, $\varepsilon=0.04$, $\sigma=0.1$ and $h=0.01$, at the different times 2,
4, 6 and 8 for $N=10000$ particles. Figure \ref{fig:alpha_moyen} gives the same simulation for~$\alpha=1.$
In the case $\alpha<1$, and especially when $\alpha$ is small, one need
to take a very small value for the diffusion coefficient in order to
have a reasonable approximation of the solution. Indeed, the
approximation depicted on the Figure \ref{fig:alpha_petit_mauvais} is the appproximation of
the solution at times 2, 4, 6 and~8 for diffusion coefficient
$\sigma=10^{-4}$. Here, we used 10000 particles killed at a distance
$\varepsilon=0.01,$ the time step being $h=0.01$.
On Figure \ref{fig:alpha_petit_bon} we show the same simulation, with
diffusion coefficient changed to $\sigma=10^{-12}$.

\bibliography{biblio}
\bibliographystyle{plain}

\end{document}